\documentclass{amsart}

\usepackage{amsmath}
\usepackage[latin1]{inputenc}
\usepackage{amsfonts}
\usepackage{amssymb}
\usepackage{graphicx,epsfig}
\usepackage{amsthm}
\usepackage{amscd}
\usepackage[all]{xy}
\usepackage[latin1]{inputenc}

\usepackage[T1]{fontenc}

\usepackage{amsmath}

\usepackage{amsfonts}

\usepackage{amssymb}

\usepackage{amsthm}

\usepackage{graphics}

\newcommand{\mk}{\medskip}

\newcommand{\ZZ}{\mathbb{Z}}

\newcommand{\CC}{\mathbb{C}}

\newcommand{\NN}{\mathbb{N}}

\newcommand{\QQ}{\mathbb{Q}}

\newcommand{\Glie}{\mathfrak{g}}

\newcommand{\Yim}{\mathcal{Y}}

\newcommand{\Hlie}{\mathfrak{h}}

\newcommand{\demo}{\noindent {\it \small Proof:}\quad}

\renewcommand{\NN}{\ensuremath{\mathbb{N}}}

\renewcommand{\CC}{\ensuremath{\mathbb{C}}}

\renewcommand{\QQ}{\ensuremath{\mathbb{Q}}}

\newcommand{\U}{\mathcal{U}}

\newtheorem{thm}{Theorem}[section]

\newtheorem{defi}[thm]{Definition}

\newtheorem{cor}[thm]{Corollary}

\newtheorem{prop}[thm]{Proposition}

\newtheorem{lem}[thm]{Lemma}

\author{David Hernandez}
\address{David Hernandez: \'Ecole Normale Sup\'erieure - DMA, 45, Rue d'Ulm F-75230 PARIS,
  Cedex 05  FRANCE}
\email{David.Hernandez@ens.fr\\ URL: http://www.dma.ens.fr/$\sim$dhernand}

\usepackage[all]{xy}
\title{Monomials of $q$ and $q,t$-characters for non simply-laced quantum affinizations}

\begin{document}

\begin{abstract} Nakajima \cite{Naa, Nab} introduced the morphism of $q,t$-characters for finite dimensional representation of simply-laced quantum affine algebras : it is a $t$-deformation of the Frenkel-Reshetikhin's morphism of q-characters (sum of monomials in infinite variables). In \cite{her02} we generalized the construction of $q,t$-characters for non simply-laced quantum affine algebras. First in this paper we prove a conjecture of \cite{her02} : the monomials of $q$ and $q,t$-characters of standard representations are the same in non simply-laced cases (the simply-laced cases were treated in \cite{Nab}) and the coefficients are non negative. In particular those $q,t$-characters can be considered as $t$-deformations of $q$-characters. In the proof we show that for quantum affine algebras of type $A, B, C$ and quantum toroidal algebras of type $A^{(1)}$ the $l$-weight spaces of fundamental representations are of dimension 1. Eventually we show and use a generalization of a result of \cite{Fre, Fre2, Naams} : for general quantum affinizations we prove that the $l$-weights of a $l$-highest weight simple module are lower than the highest $l$-weight in the sense of monomials. \end{abstract}

\maketitle

\tableofcontents

\section{Introduction}\label{intro} V.G. Drinfel'd \cite{Dri1} and M. Jimbo \cite{jim} associated, independently, to any symmetrizable Kac-Moody algebra $\Glie$ and any complex number $q\in\CC^*$ a Hopf algebra $\U_q(\Glie)$ called quantum Kac-Moody algebra.

\noindent In this paper we suppose that $q\in\CC^*$ is not a root of unity. In the case of a semi-simple Lie algebra $\Glie$ of rank $n$ (ie. with a finite Cartan matrix), the structure of the Grothendieck ring $\text{Rep}(\U_q(\Glie))$ of finite dimensional representations of the quantum finite algebra $\U_q(\Glie)$ is well understood. It is analogous to the classical case $q=1$. 

\noindent For the general case of Kac-Moody algebras the picture is less clear. The representation theory of the quantum affine algebra $\U_q(\hat{\Glie})$ is of particular interest (see \cite{Cha, Cha2}). In this case there is a crucial property of $\U_q(\hat{\Glie})$: it has two realizations, the usual Drinfel'd-Jimbo realization and a new realization (see \cite{Dri2, bec}) as a quantum affinization of a quantum finite algebra $\U_q(\Glie)$.

\noindent To study the finite dimensional representations of $\U_q(\hat{\Glie})$ Frenkel and Reshetikhin \cite{Fre} introduced $q$-characters which encode the (pseudo)-eigenvalues of some commuting elements $\phi_{i,\pm m}^{\pm}$ ($m\geq 0$) of the Cartan subalgebra $\U_q(\hat{\Hlie})\subset\U_q(\hat{\Glie})$ (see also \cite{kn}) : for $V$ a finite dimensional representations we have :
$$V=\underset{\gamma\in \CC^{I\times \ZZ}}{\bigoplus} V_{\gamma}$$ 
where for $\gamma=(\gamma_{i,\pm m}^{(\pm)})_{i\in I,m\geq 0}$, $V_{\gamma}$ is a simultaneous generalized eigenspace ($l$-weight space):
$$V_{\gamma}=\{x\in V/\exists p\in\NN,\forall i\in I,\forall m\geq 0,(\phi_{i,\pm m}^{(\pm)}-\gamma_{i,\pm m}^{(\pm)})^p.x=0\}$$
The morphism of $q$-characters is an injective ring homomorphism:
$$\chi_q:\text{Rep}(\U_q(\hat{\Glie}))\rightarrow\Yim=\ZZ[Y_{i,a}^{\pm}]_{i\in I,a\in\CC^*}$$
where $\text{Rep}(\U_q(\hat{\Glie}))$ is the Grothendieck ring of finite dimensional (type 1)-representations of $\U_q(\hat{\Glie})$ and $I=\{1,...,n\}$, and : 
$$\chi_q(V)=\underset{\gamma\in \CC^{I\times \ZZ}}{\sum}\text{dim}(V_{\gamma})m_{\gamma}$$
where $m_{\gamma}\in\Yim$ depends of $\gamma$. In particular $\text{Rep}(\U_q(\hat{\Glie}))$ is commutative and isomorphic to $\ZZ[X_{i,a}]_{i\in I,a\in\CC^*}$.

\noindent In the finite simply laced-case (type $ADE$) Nakajima \cite{Naa, Nab} introduced $t$-analogs of $q$-characters. The motivations are the study of filtrations induced on representations by (pseudo)-Jordan decompositions, the study of the decomposition in irreducible modules of tensorial products and the study of cohomologies of certain quiver varieties. The morphism of $q,t$-characters is a $\ZZ$-linear map :
$$\chi_{q,t}:\text{Rep}(\U_q(\hat{\Glie}))\rightarrow\Yim_t=\ZZ[Y_{i,a}^{\pm},t^{\pm}]_{i\in I,a\in\CC^*}$$ 
which is a deformation of $\chi_q$ and multiplicative in a certain sense. A combinatorial axiomatic definition of $q,t$-characters is given. But the existence is non-trivial and is proved with the geometric theory of quiver varieties which holds only in the simply laced case. 

\noindent In \cite{her02} we defined and constructed $q,t$-characters for a finite (non necessarily simply-laced) Cartan matrix $C$ with a new approach motivated by the non-commutative structure of $\U_q(\hat{\Hlie})\subset\U_q(\hat{\Glie})$, the study of screening currents of \cite{Freb} and of deformed screening operators $S_{i,t}$ of \cite{her01}. It generalizes the construction of Nakajima to non-simply laced cases.

\noindent The quantum affinization process (that Drinfel'd \cite{Dri2} described for constructing the second realization of a quantum affine algebra) can be extended to all symmetrizable quantum Kac-Moody algebras $\U_q(\Glie)$ (see \cite{jin, Naams, her04}). One obtains a new class of algebras called quantum affinizations : the quantum affinization of $\U_q(\Glie)$ is denoted by $\U_q(\hat{\Glie})$. It has a triangular decomposition \cite{her04}. For example the quantum affinization of a quantum affine algebra is called a quantum toroidal algebra. The quantum affine algebras are the simplest examples and are very special because they are also quantum Kac-Moody algebras. In the following, general quantum affinization means with an invertible quantum Cartan matrix (it includes most interesting cases like affine and toroidal quantum affine algebras, see section \ref{qcar}). In \cite{her04} we developed the representation theory of general quantum affinizations and constructed a generalization of the $q$-characters morphism which appears to be a powerful tool for this investigation. In particular we proved that the new Drinfel'd coproduct leads to the construction of a fusion product on the Grothendieck group.

The results of this paper can be divided in three parts :

1) First we prove that for general quantum affinizations, the $l$-weights $m'\in\Yim$ of a simple module of $l$-highest weight $m\in\Yim$ are lower than $m$ in the sense of monomials (theorem \ref{cetconj}) : it means that $m'm^{-1}$ is a product of certain $A_{i,l}^{-1}\in\Yim$. For $C$ finite, this result was conjectured and partly proved in \cite{Fre} and proved in \cite{Naams} ($ADE$-case) and \cite{Fre2} (finite case). In the general case no universal $\mathcal{R}$-matrix has been defined : so we propose a new proof based on the study of $\U_q(\hat{sl_2})$-Weyl modules introduced in \cite{chaw1}. This first result is used in the proof of the following points :

2) We prove a conjecture of \cite{her02} : let $\U_q(\hat{\Glie})$ be a quantum affine algebra ($C$ finite) and $M$ be a standard module of $\U_q(\hat{\Glie})$ (tensorial product of fundamental representations). We prove that the coefficients of $\chi_{q,t}(M)$ are in $\NN[t^{\pm}]$ and that the monomials of $\chi_{q,t}(M)$ are the monomials of $\chi_q(M)$ (theorem \ref{posqt}) (the case $ADE$ follows from the work of Nakajima \cite{Nab}; in this paper the non-simply laced case is treated.) In particular the $q,t$-characters for quantum affine algebras have a finite number of monomials : this result shows that the $q,t$-characters of \cite{her02} can be considered as $t$-deformations of $q$-characters for all quantum affine algebras. In particular it is an argument for the existence of a geometric model behind the $q,t$-characters in non simply-laced cases (in the simply laced-case the standard module can be realized in the K-theory of quivers varieties).

3) In the proof of the conjecture we study combinatorial properties of $q$-characters : we prove that for quantum affinizations of type $A, A^{(1)}, B, C$ the $l$-weight spaces of fundamental representations are of dimension 1 (theorem \ref{mainprem}). Note that this property is not true in general, for example for type $D$.

\noindent Our proof is based on an investigation of the classical algorithm (see \cite{Fre2, her03}) which gives $q$-characters. The proof is direct without explicit computation. Note that for type $A,B,C$ the result could follow from explicit computation of the specialized $\mathcal{R}$-matrix, as explained in \cite{Fre}. The result should produce the formulas of \cite{Frea}; however with this method it would not be easy to decide if the coefficients are 1 (for example it is not the case for type $D_4$). Moreover it allows us to extend the proof to quantum toroidal algebras of type $A^{(1)}$.

This paper is organized as follows : in section \ref{rmd} we give reminders on representations of quantum affinizations and their $q$-characters. In section \ref{mres} we state and prove theorem \ref{cetconj} (the $l$-weights of a $l$-highest weight simple module are lower than the highest $l$-weight in the sense of monomials) and state theorem \ref{mainprem} (on $q$-characters of fundamental representations) and give technical complements. The proof of theorem  \ref{mainprem} is based on a case by case investigation explained in sections \ref{typea}, \ref{typec}, \ref{typeb}. In section \ref{qtapp} we give reminders on $q,t$-characters and we prove theorem \ref{posqt} (on coefficients of $q,t$-characters of standard monomials). For the theorem \ref{posqt} type $F_4$, our proof is based on results obtained by a computer program written with Travis Schedler, and the results are given in the appendix (section \ref{anncalc}).

\noindent {\bf Acknowledgments} : the author would like to thank Travis Schedler for the computer program we wrote.

\section{Reminder}\label{rmd}

\subsection{Representations of quantum affinizations}

Let $C=(C_{i,j})_{1\leq i,j\leq n}$ be a symmetrizable (non necessarily finite) Cartan matrix and $I=\{1,...,n\}$. Let $D=\text{diag}(r_1,...,r_n)$ such that $B=DC$ is symmetric. We consider $(\Hlie, \Pi, \Pi^{\vee})$ a realization of $C$, the weight lattice $P\subset \Hlie^*$, the roots $\alpha_1,...,\alpha_n\in P$, the set of dominant weights $P^+\subset P$, the relation $\leq$ on $P$, the map $\nu : \Hlie^*\rightarrow \Hlie$ (see \cite{her04}). 

\noindent Let $q\in\CC^*$ not a root of unity. Let $\U_q(\Glie)$ be the quantum Kac-Moody algebra of Cartan matrix $C$. Let $\U_q(\hat{\Glie})\supset \U_q(\Glie)$ be the quantum affinization of $\U_q(\Glie)$, with generators $x_{i,r}^{\pm}, k_h, c^{\pm \frac{1}{2}}, \phi_{i,\pm m}^{(\pm)}$, where $i\in I, r\in\ZZ, m\geq 0, h\in\Hlie$ (see for example \cite{her04}). A $\U_q(\hat{\Glie})$-module is said to be integrable if it is integrable as a $\U_q(\Glie)$-module.

\noindent Denote by $P_l$ the set of $l$-weights, that is to say of couple $(\lambda,\Psi)$ such that $\lambda\in \Hlie^*$, $\Psi=(\Psi_{i,\pm m}^{\pm})_{i\in I, m\geq 0}$, $\Psi_{i,\pm m}^{\pm}\in\CC$, $\Psi_{i,0}^{\pm}=q_i^{\pm \lambda(\alpha_i^{\vee})}$. Note that if $C$ is finite, $\lambda$ is uniquely determined by $\Psi$.

\noindent A $\U_q(\hat{\Glie})$-module $V$ is said to be of $l$-highest weight $(\lambda, \Psi)\in P_l$ if there is $v\in V$ such that ($i\in I, r\in\ZZ, m\geq 0, h\in \Hlie$):
$$x_{i,r}^+.v=0\text{ , }V=\U_q(\hat{\Glie}).v\text{ , }\phi_{i,\pm m}^{\pm}.v=\Psi_{i,\pm m}^{\pm}v\text{ , }k_h.v=q^{\lambda(h)}.v$$
For $(\lambda,\Psi)\in P_l$, let $L(\lambda,\Psi)$ be the simple $\U_q(\hat{\Glie})$-module of $l$-highest weight $(\lambda,\Psi)$ (see \cite{her04}).

\noindent Let $P_l^+$ be the set of dominant $l$-weights, that is to say  the set of $(\lambda,\Psi)\in P_l$ such that there exist (Drinfel'd)-polynomials $P_i(z)\in\CC[z]$ ($i\in I$) of constant term $1$ such that in $\CC[[z]]$ (resp. in $\CC[[z^{-1}]]$):
$$\underset{m\geq 0}{\sum} \Psi_{i,\pm m}^{\pm} z^{\pm m}=q_i^{\text{deg}(P_i)}\frac{P_i(zq_i^{-1})}{P_i(zq_i)}$$

\begin{thm}\label{simpint} For $(\lambda, \Psi)\in P_l$, $L(\lambda, \Psi)$ is integrable if and only $(\lambda, \Psi)\in P_l^+$.\end{thm}

\noindent If $\Glie$ is finite (case of a quantum affine algebra) it is a result of Chari-Pressley in \cite{Cha, Cha2}. Moreover in this case the integrable $L(\lambda,\Psi)$ are finite dimensional. If $C$ is simply-laced the result is proved by Nakajima in \cite{Naams}. If $C$ is of type $A_n^{(1)}$ (toroidal $\hat{sl_n}$-case) the result is proved by Miki in \cite{mi}. In general the result is proved in \cite{her04}.

\noindent Denote by $\text{Rep}(\U_q(\hat{\Glie}))$ the Grothendieck group of decomposable integrable representations of type 1 (see \cite{her04}). The operators $k_h, \phi^{(\pm)}_{i,\pm m}\in \U_q(\hat{\Glie})$ ($h\in\Hlie,i\in I, m\in\ZZ$) commute on $V\in\text{Rep}(\U_q(\hat{\Glie}))$. So we have a $l$-weight space decomposition: 
$$V=\underset{(\lambda,\gamma)\in P_l}{\bigoplus} V_{\lambda, \gamma}$$ 
$$V_{\lambda, \gamma}=\{x\in V_{\lambda}/\exists p\in\NN,\forall i\in I,\forall m\geq 0,(\phi_{i,\pm m}^{(\pm)}-\gamma_{i,\pm m}^{(\pm)})^p.x=0\}\subset V_{\lambda}=\{v\in V/ \forall h\in\Hlie, k_h.v=q^{\lambda(h)}v\}$$

\noindent Let $QP_l^+\subset P_l$ be the set of $(\mu,\gamma)\in P_l$ such that there exist polynomials $Q_i(z), R_i(z)\in\CC[z]$ ($i\in I$) of constant term $1$ such that in $\CC[[z]]$ (resp. in $\CC[[z^{-1}]]$):
$$\underset{m\geq 0}{\sum}\gamma_{i,\pm m}^{\pm} z^{\pm m}= q_i^{\text{deg}(Q_i)-\text{deg}(R_i)}\frac{Q_i(zq_i^{-1})R_i(zq_i)}{Q_i(zq_i)R_i(zq_i^{-1})}$$

\noindent In particular $P_l^+\subset QP_l^+$.

\begin{prop}\label{fr} Let $V$ be a module in $\text{Rep}(\U_q(\hat{\Glie}))$ and $(\mu,\gamma)\in P_l$. If $\text{dim}(V_{\mu, \gamma})>0$ then $(\mu,\gamma)\in QP_l^+$.\end{prop}

\noindent The result is proved in \cite{Fre} for $C$ finite. The generalization is straightforward (see \cite{her04}).

\subsection{q-characters}\label{qcar} Let $z$\label{z} be an indeterminate. We denote $z_i=z^{r_i}$ and for $l\in\ZZ$, $[l]_z=\frac{z^l-z^{-l}}{z-z^{-1}}\in\ZZ[z^{\pm}]$. Let $C(z)$ be the quantized Cartan matrix defined by ($i\neq j\in I$):
$$C_{i,i}(z)=z_i+z_i^{-1}\text{ , }C_{i,j}(z)=[C_{i,j}]_z$$
In the rest of this paper we suppose that $C(z)$ is invertible, that is to say $\text{det}(C(z))\neq 0$. We have seen in lemma 6.9 of \cite{her03} that the condition $(C_{i,j}<-1\Rightarrow -C_{j,i}\leq r_i)$ implies that $\text{det}(C(z))\neq 0$. In particular finite and affine Cartan matrices (where we impose $r_1=r_2=2$ for $A_1^{(1)}$) satisfy this condition. 

\noindent Consider formal variables $Y_{i,a}^{\pm}$ ($i\in I, a\in\CC^*$) and $k_{\omega}$ ($\omega\in\Hlie$) ($k_0=1$). Let $A$ be the commutative group of monomials of the form $m=\underset{i\in I,a\in\CC^*}{\prod}Y_{i,a}^{u_{i,a}(m)} k_{\omega(m)}$ where a finite number of $u_{i,a}(m)\in\ZZ$ are non zero, $\omega(m)\in \Hlie$ (the coweight of $m$), and such that for $i\in I$ : $\alpha_i(\omega(m))=r_i\underset{a\in\CC^*}{\sum}u_{i,a}(m)$.

\noindent For $(\mu, \Gamma)\in QP_l^+$ we define $Y_{\mu,\Gamma}\in A$ by:
$$Y_{\mu, \Gamma}=k_{\nu(\mu)}\underset{i\in I, a\in\CC^*}{\prod}Y_{i,a}^{\beta_{i,a} -\gamma_{i,a}}$$
where $\beta_{i,a}, \gamma_{i,a}\in\ZZ$ are defined by $Q_i(u)=\underset{a\in\CC^*}{\prod}(1-ua)^{\beta_{i, a}}$ , $R_i(u)=\underset{a\in\CC^*}{\prod}(1-ua)^{\gamma_{i, a}}$. 

\noindent For $(\mu,\Gamma)\in QP_L^+$ and $V$ a $\U_q(\hat{\Glie})$-module, we denote $V_m=V_{\mu,\Gamma}$ where $m=Y_{\mu,\Gamma}$.

\noindent For $\chi\in A^{\ZZ}$ we say $\chi\in\Yim$ if for $\lambda\in\Hlie$, there is a finite number of monomials of $\chi$ such that $\omega(m)=\lambda$ and there is a finite number of element $\lambda_1,...,\lambda_s\in \Hlie^*$ such that the coweights of monomials of $\chi$ are in $\underset{j=1...s}{\bigcup}\nu(D(\lambda_j))$ (where $D(\lambda_j)=\{\mu\in\Hlie^*/\mu\leq \lambda_j\}$). In particular $\Yim$ has a structure of $\Hlie$-graded $\ZZ$-algebra.

\begin{defi} For $V\in\text{Rep}(\U_q(\hat{\Glie}))$ a representation, the $q$-character\label{chiqdefi} of $V$ is:
$$\chi_q(V)=\underset{(\mu,\Gamma)\in QP_l^+}{\sum} \text{dim}(V_{\mu,\Gamma})Y_{\mu,\Gamma}\in\Yim$$\end{defi}

\noindent If $C$ is finite the construction is given in \cite{Fre} and it is proved that $\chi_q:\text{Rep}(\U_q(\hat{\Glie}))\rightarrow \Yim$ is an injective ring homomorphism (with the ring structure on $\text{Rep}(\U_q(\hat{\Glie}))$ deduced from the Hopf algebra structure of $\U_q(\hat{\Glie})$). 

\noindent For general $C$, $\chi_q$ is defined in \cite{her04}. A priori there is no ring structure on $\Yim$ that comes from a tensor product, but we proved \cite{her04}:

\begin{thm} The image $\text{Im}(\chi_q)\subset\Yim$ is a subalgebra of $\Yim$.\end{thm}

\noindent Let $*$ be the induced commutative product on $\text{Im}(\chi_q)\subset\Yim$. Using a deformation of the new Drinfel'd coproduct we proved in \cite{her04} :

\begin{thm}\label{h4thm} For $(\lambda,\Psi), (\lambda',\Psi')\in P_l^+$ we have :
$$L(\lambda,\Psi)*L(\lambda',\Psi')=L(\lambda+\lambda',\Psi\Psi')+\underset{(\mu,\Phi)\in P_l^+/ \mu < \lambda +\lambda'}{\sum}Q_{\lambda,\Psi,\lambda',\Psi'}(\mu,\Phi)L(\mu,\Phi)$$
where the integers $Q_{\lambda,\Psi,\lambda',\Psi'}(\mu,\Phi)\geq 0$.\end{thm}

\subsection{Notations and technical tools}\label{notst} For $i\in I$ and $a\in\CC^*$ we set:
$$A_{i,a}=k_{r_i\alpha_i^{\vee}}Y_{i,aq_i^{-1}}Y_{i,aq_i}\underset{j/C_{j,i}<0}{\prod}\underset{r=C_{j,i}+1,C_{j,i}+3,...,-C_{j,i}-1}{\prod}Y_{j,aq^r}^{-1}\in A$$
The $A_{i,l}$ are algebraically independent (see \cite{her02}). Let $\mathcal{A}=\ZZ[A_{i,a}^{-1}]_{i\in I,a\in\CC^*}\subset\Yim$. 

\noindent For a product $M\in A$ such that $\omega(M)\in\ZZ\alpha_1\oplus ...\oplus\ZZ\alpha_n$, denote $\omega(M)=-v_1(M)\alpha_1-...-v_n(M)\alpha_n$ and $v(M)=v_1(M)+...+v_n(M)$. $v$ defines a $\NN$-gradation on $\mathcal{A}$. 

\begin{defi}\label{ordrmonom} For $m,m'\in A$, we say that $m\geq m'$ if $m'm^{-1}\in\mathcal{A}$.\end{defi}

For $m\in A$ and $J\subset I$, denote $u_J(m)=\underset{j\in J, a\in\CC^*}{\sum}u_{j,a}(m)$, $m^{(J)}=k_{\omega (m)}\underset{j\in J, a\in\CC^*}{\prod}Y_{j,a}^{u_{j,a}(m)}$ and ($j\in I$):
$$u_j^{\pm}(m)=\pm \underset{l\in\ZZ/ \pm u_{j,l}(m)>0}{\sum}u_{j,l}(m)\text{ , }u_J^{\pm}(m)=\underset{j\in J}{\sum}u_j^{\pm}(m)$$
For $J\subset I$, denote $B_J\subset A$ the set of $J$-dominant monomials (ie $\forall j\in J, l\in\ZZ$, $u_{j,l}(m)\geq 0$) and $B=B_I$. Note that for $(\lambda,\Psi)\in QP_l^+$ : ($(\lambda,\Psi)\in P_l^+\Leftrightarrow Y_{\lambda,\Psi}\in B$).

\noindent For $m\in B$ denote $V_m=L(\lambda,\Psi)\in\text{Rep}(\U_q(\hat{\Glie}))$ where $(\lambda,\Psi)\in P_l^+$ is given by $Y_{\lambda,\Psi}=m$. In particular for $i\in I,a\in\CC^*$ denote $V_i(a)=V_{k_{\nu(\Lambda_i)}Y_{i,a}}$ and $X_{i,a}=\chi_q(V_{i,a})$. The simple modules $V_i(a)$ are called fundamental representations. 

\noindent Denote $M_m=\underset{i\in I,a\in\CC^*}{\prod}V_{i,a}^{* u_{i,a}(m)}\in\text{Rep}(\U_q(\hat{\Glie}))$. We have $\chi_q(M_m)=\underset{i\in I,a\in\CC^*}{\prod}X_{i,a}^{u_{i,a}(m)}$.

For $J\subset I$ we denote by $\Glie_J$ the Kac-Moody algebra of Cartan matrix $(C_{i,j})_{i,j\in J}$ and by $\chi_q^J$ the morphism of $q$-characters for $\U_q(\hat{\Glie}_J)\subset\U_q(\hat{\Glie})$. Let us recall the definition of the morphism $\tau_J$ (section 3.3 in \cite{Fre2} for finite case and \cite{her04} for general case) :

\noindent We suppose that $\Glie_J$ is finite. Let $\Hlie_J^{\perp}=\{\omega\in\Hlie/\forall i\in J,\alpha_i(\omega)=0\}$ and $\Hlie_J=\underset{i\in J}{\bigoplus}\QQ \Lambda_i^{\vee}$. Consider formal variables $k'_{\omega}$ ($\omega\in\Hlie_J$), $k_{\omega}$ ($\omega\in\Hlie_J^{\perp}$), $Y_{i,a}^{\pm}$ ($i\in J, a\in\CC^*$), $Z_{j,c}$ ($j\in I-J$, $c\in\CC^*$). Let $A^{(J)}$ be the commutative group of monomials :
$$m=k'_{\omega'(m)}k_{\omega(m)}\underset{i\in J,a\in\CC^*}{\prod}Y_{i,a}^{u_{i,a}(m)}\underset{j\in I-J,c\in\CC^*}{\prod}Z_{j,c}^{z_{j,c}(m)}$$ 
where a finite number of $u_{i,a}(m),z_{j,c}(m),r(m)\in\ZZ$ are non zero, $\omega(m)\in \Hlie_J^{\perp}$ and such that for $i\in J$, $\alpha_i(\omega'(m))=r_iu_i(m)=r_i\underset{a\in\CC^*}{\sum}u_{i,a}(m)$.

\noindent Let $\tau_J:A\rightarrow A^{(J)}$ be the group morphism defined formally by ($j\in I$, $a\in\CC^*$, $\lambda\in \Hlie$):
$$\tau_J(Y_{j,a})=Y_{j,a}^{\epsilon_{j,J}}\underset{k\in I-J}{\prod}\underset{r\in\ZZ}{\prod}Z_{k,a q^r}^{p_{j,k}(r)}\text{ , }\tau_J(k_{\lambda})=k'_{\underset{i\in J}{\sum}\alpha_i(\lambda)\Lambda_i^{\vee}}k_{\lambda-\underset{i\in J}{\sum}\alpha_i(\lambda)\Lambda_i^{\vee}}$$
where $j\in J\Leftrightarrow\epsilon_{j,J}=1$ and $j\notin J\Leftrightarrow\epsilon_{j,J}=0$. The $p_{i,j}(r)\in\ZZ$ are defined in the following way : we write $\tilde{C}(z)=\frac{\tilde{C}'(z)}{d(z)}$ where $d(z), \tilde{C}'_{i,j}(z)\in\ZZ[z^{\pm}]$ and $(D(z)\tilde{C}'(z))_{i,j}=\underset{r\in\ZZ}{\sum}p_{i,j}(r)z^r$.

\noindent It is proved in \cite{Fre2} (finite case) and in \cite{her04} (the proof is given for the $\tau_{\{i\}}$ ($i\in I$), but the proof for $\tau_J$ ($J\subset I$, $\Glie_J$ finite) is the same) :

\begin{lem}\label{aidedeux} Consider $V$ a module in $\text{Rep}(\U_q(\hat{\Glie}))$ and a decomposition $\tau_J(\chi_q(V))=\underset{k}{\sum}P_kQ_k$ where  $P_k\in\ZZ[Y_{i,a}^{\pm},k_{h}']_{i\in J, a\in\CC^*, h\in\Hlie_J}$, $Q_k$ is a monomial in $\ZZ[Z_{j,c}^{\pm}, k_h]_{j\in I-J,c\in\CC^*,h\in \Hlie_J^{\perp}}$ and all monomials $Q_k$ are distinct. Then the restriction of $V$ to $\U_q(\hat{\Glie}_J)$ is isomorphic to $\underset{k}{\bigoplus}V_k$ where $V_k$ is a $\U_q(\hat{\Glie}_J)$-module and $\chi_q^J(V_k)=P_k$.\end{lem}

\subsection{Classical algorithm}\label{classalgo} Consider $\mathfrak{K}=\underset{i\in I}{\bigcap}\mathfrak{K}_i\subset\Yim$ where $\mathfrak{K}_i=\text{Ker}(S_i)\subset\Yim$ is the kernel of the screening operator $S_i$ (see \cite{her04}). 

\begin{thm}\label{kim} We have $\mathfrak{K}=\text{Im}(\chi_q)$ and it is a subalgebra of $\Yim$.\end{thm}

\noindent The result in proved in \cite{Fre2} for $C$ finite and in \cite{her04} in general. Note that for $m\in B_i$, there is a unique $F_i(m)\in\mathfrak{K}_i$ such that $m$ is the unique $i$-dominant monomial of $F_i(m)$ (see \cite{her02}).

\noindent  In \cite{her02} a classical algorithm (and also a $t$-deformation of it) is proposed : if it is well-defined, it gives for $m\in B$ a $F(m)\in\mathfrak{K}$ such that $m$ is the unique dominant monomial of $F(m)$. Such an algorithm was first used in \cite{Fre2} for finite case (see also \cite{her03}). Note that if $F(m)$ exists, it is unique (see \cite{her02}). Let us describe this algorithm : first for $m\in B$ we have to define the set $D_m$ :

\begin{defi} For $m\in B$, we say that $m'\in D_m$ if and only if there is a finite sequence $(m_0=m,m_1,...,m_R=m')$, such that for all $1\leq r\leq R$, there is $j\in I$ such that $m_{r-1}\in B_j$ and $m_r$ is a monomial of $F_j(m_{r-1})$.\end{defi}

\noindent In particular the set $D_m$ is countable (see \cite{her02}) and $m'\in D_m\Rightarrow m'\leq m$. Denote $D_m=\{m_0=m, m_1,m_2,...\}$ such that $m_r\leq m_{r'}\Rightarrow r\geq r'$.

\noindent For $r,r'\geq 0$ and $j\in I$ denote $[F_j(m_{r'})]_{m_r}\in\ZZ$ the coefficient of $m_r$ in $F_j(m_{r'})$. 

\noindent We call classical algorithm the following inductive definition of the sequences $(s(m_r))_{r\geq 0}\in\ZZ^{\NN}$, 
\\$(s_j(m_r))_{r\geq 0}\in\ZZ^{\NN}$ ($j\in I$) : $s(m_0)=1$ , $s_j(m_0)=0$ and for $r\geq 1, j\in I$:
$$s_j(m_r)=\underset{r'<r}{\sum}(s(m_{r'})-s_j(m_{r'}))[F_j(m_{r'})]_{m_r} $$
$$m_r\notin B_j\Rightarrow s(m_r)=s_j(m_r)\text{ , }m_r\in B\Rightarrow s(m_r)=0$$
It follows from theorem \ref{kim} that the classical algorithm is weel-defined and for all $m\in B$, $F(m)\in\mathfrak{K}$ exists (see section 5.5.4 in \cite{her04}).

\section{Monomials of $q$-characters}\label{mres} In this section we state the two main results on $q$-characters of this paper :  theorems \ref{cetconj} and \ref{mainprem}.

\subsection{First result} In this section we prove that for $m'$ a $l$-weight of $V_m$ we have $m'\leq m$ (theorem \ref{cetconj}). This result is proved in \cite{Fre, Fre2} for $C$ finite. In the general case a universal $\mathcal{R}$-matrix has not been defined so we propose a new proof based on the Weyl modules introduced in \cite{chaw1}.

\begin{defi}\label{dun} For $m\in B$, denote $L(m)=\chi_q(V_m)$ and by $D(m)$ the set of monomials of $L(m)$.\end{defi}

\noindent The partial order on monomials is set in definition \ref{ordrmonom}.

\begin{thm}\label{cetconj} For $m\in B$ and $m'\in D(m)$, we have $m'\leq m$.\end{thm}

\noindent In this section \ref{mres} we prove this theorem. First let us show some lemmas which will be useful :

\begin{lem}\label{subh} Let $V$ be a $\U_q(\hat{\Glie})$-module and $W\subset V$ a $\U_q(\hat{\Hlie})$-submodule of $V$. Then for $i\in I$, $W_i'=\underset{r\in\ZZ}{\sum}x_{i,r}^-.W$ is a $\U_q(\hat{\Hlie})$-submodule of $V$.\end{lem}

\demo For $w\in W$, $j\in J$, $m,r\in\ZZ$ ($m\neq 0$), $h\in\Hlie$ we have :
$$h_{j,m}.(x_{i,r}^-.w)=x_{i,r}^-.(h_{j,m}.w)-\frac{1}{m}[mB_{i,j}]_q x_{i,m+r}^-.w\in W_i'$$
$$k_h.(x_{i,r}^-.w)= x_{i,r}^-.(q^{\alpha_i(h)}k_h.w)\in W_i'$$
\qed

\noindent Note that $q$-character of an (integrable) $\U_q(\hat{\Hlie})$-module is well-defined (see section 5.4 of \cite{her04}).

\begin{lem}\label{cassl2} Suppose that $\Glie=sl_2$ and let $L$ be a finite dimensional $\U_q(\hat{\Glie})$-module ($\Lambda^{\vee}$ is the fundamental coweight).

(i) If $L$ is of $l$-highest weight $M$ then $L_{m'}\neq \{0\}$ implies $m'\leq M$.

(ii) For $p\in\ZZ$, let $L_p=\underset{\lambda\in P^*/\lambda(\Lambda^{\vee})\geq p}{\sum}L_{\lambda}$ and $L'_p=\underset{r\in\ZZ}{\sum}x_r^-.L_p$. Then $L_p, L_p'$ are $\U_q(\hat{\Hlie})$-submodule of $L$ and $(L_p')_{m'}\neq 0\Rightarrow \exists m, m'\leq m$ and $(L_p)_{m}\neq \{0\}$.\end{lem}

\demo (i) Consider the Weyl module $W_q(M)$ of $l$-highest weight $M$ defined in \cite{chaw1} : $W_q(M)$ is the universal finite dimensional module of $l$-highest weight $M$ such that all finite dimensional module of highest $l$-weight $M$ is a quotient of $W_q(M)$. In particular $L$ is a quotient of $W_q(M)$. So it suffices to study $W_q(M)$. For $\U_q(\hat{sl_2})$, the Weyl modules are explicitly described in \cite{chaw2} : in particular the dimension of $W_q(M)$ is $2^m$ where $m=u(M)=\underset{i\in I, a\in\CC^*}{\sum}u_{i,a}(M)$. But (see \cite{var2, Aka, Fre2}) there is a standard module (tensorial product of fundamental representation) of highest $l$-weight $M$. The dimension of such a standard module is $2^m$ and it is a quotient of $W_q(M)$. So $W_q(M)$ is isomorphic to a standard module. The $q$-character of a standard module is known (see section \ref{notst}), in particular for a $l$-weight $m'$ of $W_q(M)$ we have $m'\leq M$.

(ii) $L_p$ is a $\U_q(\hat{\Hlie})$-submodule of $L$ because the action of $\U_q(\hat{\Hlie})$ does not change the weight, so it follows from lemma \ref{subh} that $L_p'$ is a $\U_q(\hat{\Hlie})$-submodule of $L$. Let us prove the second point by induction on $\text{dim}(L_p)$ : if $L_p=\{0\}$ we have $L_p'=\{0\}$. In general let $v$ be a $l$-highest weight vector of $L_p$ (there is at least one, see the proof of proposition 5.2 in \cite{her04}) and denote by $M$ his $l$-weight. Consider $V=\U_q(\hat{\Glie}).v$. It is a $l$-highest weight module and so it follows from (i) that $V_m\neq\{0\}\Rightarrow m\leq M$. We can use the induction hypothesis with $L^{(1)}=L/V$ and we get the result because $\chi_q(L)=\chi_q(V)+\chi_q(L^{(1)})$.\qed

End of the proof of theorem \ref{cetconj} :

\noindent We prove the result by induction on $v(m'm^{-1})\geq 0$. For $v(m'm^{-1})=0$ we have $m'=m$. In general suppose that the result is known for $v(m'm^{-1})\leq p$ and let $W=\underset{m'/v(m'm^{-1})\leq p}{\sum}(V_m)_{m'}$. Note that $W$ is a $\U_q(\hat{\Hlie})$-submodule of $V_m$. It follows from the triangular decomposition of $\U_q(\hat{\Glie})$ (see \cite{her04}) that : 
$$\underset{m'/v(m'm^{-1})=p+1}{\bigoplus}(V_m)_{m'}\subset \underset{i\in I}{\sum}W_i'\text{ where }W_i'=\underset{r\in\ZZ}{\sum}x_{i,r}^-.W$$
\noindent For $i\in I$, $W_i'$ is a $\U_q(\hat{\Hlie})$-submodule of $V_m$ (lemma \ref{subh}). In particular $W_i'=\underset{m'}{\bigoplus}(W_i'\cap (V_m)_m')=\underset{m'}{\bigoplus}(W_i)_{m'}$ and it suffices to show that for $i\in I$, $(W_i')_{m'}\neq \{0\}\Rightarrow m'\leq m$.
 
\noindent Consider the decomposition of lemma \ref{aidedeux} with $J=\{i\}$: $V=\underset{k}{\bigoplus}V_k$. We have $W=\underset{k}{\bigoplus}(V_k\cap W)$ and so $W_i'=\underset{k}{\bigoplus}(V_k\cap W_i')$ (because $V_k$ is a sub $\U_q(\hat{\Glie}_i)$-module of $V_m$). So we can use the (ii) of lemma \ref{cassl2} to the $\U_{q_i}(\hat{sl_2})$-module $V_k$ with $p_k$ such that $(V_k)_{p_k}=V_k\cap W$. We get that for $m$ a monomial of $\chi_q^i(W_i')$ there are $m'$ a monomial $\chi_q^i(W)$ and $m''\in\ZZ[Y_{i,a}^{-1}Y_{i,aq_i^2}^{-1}(k_{2r_i}^{(i)})^{-1}]_{a\in\CC^*}$ such that $m=m'm''$ (the $k_{r}^{(i)}$ are the $k_h'$ for $\U_{q_i}(\hat{sl_2})$, see \cite{her04}). It follows from the lemma 5.9 of \cite{her04} (see also \cite{Fre2}) that $\tau_i(A_{i,aq_i})=Y_{i,a}Y_{i,aq_i^2}k_{2r_i}^{(i)}$. So for $M$ a monomial of $\chi_q(W_i')$ there is $M'$ a monomial of $\chi_q(W)$ such that $M\leq M'$.\qed

\subsection{Second result}

\begin{thm}\label{mainprem} Let $\Glie$ be of type $A_n$ ($n\geq 1$), $A_l^{(1)}$ ($l\geq 2$), $B_n$ ($n\geq 2$) or  $C_n$ ($n\geq 2$). Let $i\in I, a\in\CC^*$. Then for $m\in D(Y_{i,a})$, for all $j\in I, l\in\ZZ$, $u_{j,l}(m)\leq 1$. In particular all coefficients of $L(Y_{i,a})$ are equal to $1$ and all $l$-weight space of $V_i(a)$ are of dimension 1.\end{thm}

\noindent The last part of the result for type $A_n$ is established in \cite{Nac}. 

\noindent Note that for type $D_n$ the statement is false : for example for the type $D_4$, the monomial $Y_{2,2}Y_{2,4}^{-1}$ has a coefficient $2$ in $\chi_q(V_2(q^0))$ (see the figure 1 in \cite{Naa}). For type $F_4$ it is also false (see section \ref{anncalc}).

\noindent Let us explain the main points of the proof : it is based on the study of the classical algorithm in a case by case investigation : for type $A_n$ a proof is given in \cite{her02} and recalled in section \ref{typea}. The result for type $A_l^{(1)}$ is proved in section \ref{typea}, the result for type $B_n$ is proved in section \ref{typec}, the result for type $C_n$ in section \ref{typeb}. In each case we suppose the existence of a $m\in D(Y_{i,a})$, such that there is $j\in I, l\in\ZZ$, $u_{j,l}(m)\geq 2$. The classical algorithm starts from the highest weight monomial. In our proof we look at a monomial $m$ with $u_{j,l}(m)\geq 2$ and show that inductively that it implies a condition on some monomials of higher weight. In particular it leads to a contradiction on the highest weight monomial.

\noindent Note that for type $A_n$, $B_n$, $C_n$ the result could follow from explicit computation of $\chi_q(V_i(a))$. We would have to compute the specialized $\mathcal{R}$-matrix, as explained in \cite{Fre}. The result should produce the formulas of \cite{Frea}. However with this method it would not easy to decide if the coefficients are 1 (for example it is not the case for type $D_4$). In this paper the proof is direct without explicit computation. In particular it allows us to extend the proof to $A_l^{(1)}$.

\subsection{Notations}\label{notadeux} In the following (sections \ref{mres}, \ref{typea}, \ref{typec}, \ref{typeb}) we can forget the terms $k_{\lambda}$ because we work in a set $D(m)$ or $D_m$ : indeed $m'$ such that $m'\leq m$ is uniquely determined by $m$ and the $v_{i,l}(m'm^{-1})$.

\noindent For $J\subset I$, $j\in J, a\in\CC^*$ consider $A_{j,a}^{J, \pm}=(A_{j,a}^{\pm})^{(J)}$. Define $\mu_J^I:\ZZ[A_{j,a}^{J, \pm}]_{j\in J,a\in\CC^*}\rightarrow \ZZ[A_{j,a}^{\pm}]_{j\in J,a\in\CC^*}$ the ring morphism such that $\mu_J^I(A_{j,a}^{J, \pm})=A_{j,a}^{\pm}$. For $m\in B_J$, denote $L^J(m^{(J)})$ defined for $\Glie_J$ ($\Glie_J$ is the Kac-Moody algebra of Cartan matrix $(C_{i_1,i_2})_{i_1,i_2\in J}$). Define :
$$L_J(m)=m^{(I-J)}\mu_J^I((m^{(J)})^{-1}L^J(m^{(J)}))$$

\begin{defi} For $J\subset I$ and $m\in B_J$, denote by $D_J(m)$ the set of monomials of $L_J(m)$.\end{defi}

\noindent For $J=\{i\}$ and $m\in B_i$, an explicit description of $D_i(m)$ is given in \cite{Fre} : a $\sigma\subset\ZZ$ is called a $2$-segment if $\sigma$ is of the form $\sigma=\{l,l+2,...,l+2k\}$ where $l\in\ZZ, k\geq 0$. Two $2$-segment are said to be in special position if their union is a $2$-segment that properly contains each of them. All finite subset of $\ZZ$ with multiplicity $(l,u_l)_{l\in\ZZ}$ ($u_l\geq 0$) can be broken in a unique way into a union of $2$-segments which are not in pairwise special position. For $m\in B_i$ and $r\in\{1,...,2r_i\}$, consider $(\sigma_j^{(r)})_j$ the decomposition of the $(l,u_{r+2r_il}(m))_{l\in\ZZ}$ as above. Let $m^{(i)}=\underset{r=1,...,2r_i}{\prod}\underset{j}{\prod}m_{r,j}$ where $m_{r,j}=\underset{l\in \sigma_j^{(r)}}{\prod}Y_{i,r+2r_il}$, and we have :
$$D_i(m)=m^{(I-\{i\})}\underset{r=1,...,2r_i}{\prod}\underset{j}{\prod}D_i(m_{r,j})$$
where for $m=\underset{k=1...r}{\prod}Y_{i,l+2r_ik}$ :
$$D_i(m)=\{mA_{i,l+2r_ik+r_i}^{-1}, mA_{i,l+2r_ik+r_i}^{-1}A_{i,l+2r_i(k-1)+r_i}^{-1},...,mA_{i,l+2r_ik+r_i}^{-1}A_{i,l+2r_i(k-1)+r_i}^{-1}...A_{i,l+r_i}^{-1}\}$$
In particular :
\begin{lem}\label{domiseul} For $m\in B_i$ such that $\forall l\in\ZZ, u_{i,l}(m)\leq 1$, we have $F_i(m)=L_i(m)$.\end{lem}

\begin{defi}\label{fleche} Let $J\subset I$ and $i\in I, a\in\CC^*$. For $m,m'\in D(Y_{i,a})$, we denote :

$m \rightarrow_J m'$ (or $m' \leftarrow_J m$) if  $m\in B_J$ and $m'\in D_J(m)$.

$m\rightharpoonup_J m'$ (or $m'\leftharpoonup_J m$) if $v(m'Y_{i,a}^{-1})\geq v(mY_{i,a}^{-1})$ and $\exists m''\in D(Y_{i,a})$ such that $m''\rightarrow_J m$ and $m''\rightarrow_J m'$.\end{defi}

\noindent In particular $m\rightarrow_J m'$ implies $m\rightharpoonup_J m'$. For $J=\{j\}$ (resp. $J=I$) we denote $\rightarrow_j, \rightharpoonup_J$ (resp. $\rightarrow, \rightharpoonup$).

\subsection{Technical complements}

\begin{prop}\label{jdecomp} For $m\in B$ and $J\subset I$ such that $\Glie_J$ is finite, there is a unique decomposition:
$$L(m)=\underset{m'\in B_J\cap D(m)}{\sum}\lambda_J(m')L_J(m')$$
where $\lambda_J(m')\geq 0$.\end{prop}

\demo  Consider the decomposition of lemma \ref{aidedeux} with $J$ : $V=\underset{k}{\bigoplus}V_k$. We can decompose each $V_k$ in a sum of simple $\U_q(\hat{\Glie}_J)$-modules in the Grothendieck group : $\chi_q^J(V_k)=\underset{k'}{\sum}\lambda_{k,k'}L^J(m_{k,k'})$ where $m_{k.k'}\in B_J$ and $\lambda_{k,k'}\geq 0$. In particular $\tau_J^{-1}(P_kQ_k)=\underset{k'}{\sum}\lambda_{k,k'}L_J(\tau_J^{-1}(m_{k,k'}Q_k))$ (consequence of lemma 5.9 of \cite{her04}). For the uniqueness the $L_J(m')$ ($m'\in B_J$) are linearly independent.\qed

\noindent We say that a monomial $m\in B$ is right (resp. left) negative if : for $b\in\CC^*$ such that ($\exists j\in I$, $u_{j,b}(m')\neq 0$ and $\forall k\in I, l>0$ (resp. $l<0$), $u_{k,bq^l}(m')=0$), we have $\forall k\in I$, $u_{k,b}(m')\leq 0$ (see \cite{Fre2}). A product of right (resp. left) negative monomials is right (resp. left) negative.

\begin{cor}\label{propgene} For $i\in I, a\in\CC^*$ and $m'\in D(Y_{i,a})$, we have :

1) for $J\subset I$ such that $\Glie_J$ is finite, there is $m''\rightarrow_J m'$.

2) there is a finite sequence $Y_{i,a}=m_0>m_1>m_2>...>M_R=m'$ such that for all $1\leq r\leq R$, $\exists j_r\in I$, $m_{r-1}\rightarrow_{j_r}m_r$.

3) if $m'\neq Y_{i,a}$, then $m'$ is right negative

4) for $b\in\ZZ$ and $j\in I$, we have $u_{j,b}(m')\neq 0\Rightarrow b\in aq^{\ZZ}$.\end{cor}

\noindent Note that the (1) will be used intensively in the following. For $C$ finite those results are proved in \cite{Fre2}.

\demo 1) Consequence of proposition \ref{jdecomp}.

2) We use (1,3) inductively.

3) For $m'\in D_{Y_{i,a}}-\{Y_{i,a}\}$, we have $m'< Y_{i,a}$ (theorem \ref{cetconj}) and $m'$ is right or left negative, so not dominant. So as in \cite{Fre2} $m'$ is right negative.

4) As for $m\in B_j\cap \ZZ[Y_{i,aq^m}]_{m\in\ZZ}$ implies $L_j(m)\in \ZZ[Y_{i,aq^m}]_{m\in\ZZ}$ (see section \ref{notadeux}), we have $M\in B\cap \ZZ[Y_{i,aq^m}]_{m\in\ZZ}$ implies $D_m\subset \ZZ[Y_{i,aq^m}]_{m\in\ZZ}$ (see also \cite{Fre2}).\qed

\noindent As a right negative monomial is not dominant, we have :

\begin{cor} For $i\in I, a\in\CC^*$, $L(Y_{i,a})=F(Y_{i,a})$ has a unique dominant monomial $Y_{i,a}$.\end{cor}

\noindent For $c\in\CC^*$, let $\beta_c: \Yim\rightarrow\Yim$ be the ring morphism such that $\beta_c(Y_{i,a})=Y_{i,ac}$.

\begin{prop}\label{aidafm} For $a,b\in\CC^*$, $L(Y_{i,a})=\beta_{ab^{-1}}(L(Y_{i,b}))$.\end{prop}

\demo For $c\in\CC^*$, we have $\beta_c(\mathfrak{K})=\mathfrak{K}$ (see \cite{Fre2, her03}).\qed

\noindent It suffices to study (see (4) of corollary \ref{propgene} and \cite{her03}) : 
$$\chi_q:\ZZ[X_{i,q^l}]_{i\in I, l\in\ZZ}\rightarrow\ZZ[Y_{i,q^l}]_{i\in I, l\in\ZZ}$$ 
In the following we denote $\text{Rep}=\ZZ[X_{i,q^l}]_{i\in I, l\in\ZZ}$, $X_{i,l}=X_{i,q^l}$, $\Yim=\ZZ[Y_{i,q^l}^{\pm}]_{i\in I,l\geq 0}$, $Y_{i,l}^{\pm}=Y_{i,q^l}^{\pm}$. A $\text{Rep}$-monomials is a product of the $X_{i,l}$.

\begin{lem}\label{prodstan} For $m\in B$, we have : 
$$D(m)\subset\underset{j\in I,l\in\ZZ}{\prod} D(Y_{j,l})^{u_{j,l}(m)}$$\end{lem}

\demo $\underset{j\in I , l\in\ZZ}{\prod} D(Y_{j,l})^{u_{j,l}(m)}$ is the set of monomials of $\chi_q(M_m)$. Then see theorem \ref{h4thm}.\qed

\begin{lem}\label{inter} Let $m_1,m_2\in B_i$ such that $\forall l\in\ZZ$, $u_{i,l}(m_1)\leq 1$ and $u_{i,l}(m_2)\leq 1$. Then $D_i(m_1)=D_i(m_1)$ (resp. $D_i(m_2)=D_i(m_2)$) and $D_i(m_1)\cap D_i(m_2)= \emptyset\Leftrightarrow m_1\neq m_2$.\end{lem}

\demo Let us write $m_1^{(i)}=\underset{r=1,...,2r_i}{\prod}\underset{j}{\prod}m_{\sigma_j^{(r)}}^{(r)}$ as in section \ref{notadeux} . Denote $\overline{\sigma_j^{(r)}}=\sigma_j^{(r)}\cup\{\text{max}(\sigma_j^{(r)})+2r_i\}$. It follows from the hypothesis of the lemma that $(j,r)\neq (j',r')\Rightarrow \overline{\sigma_j^{(r)}}\cap\overline{\sigma_{j'}^{(r')}}=\emptyset$. Moreover for $m'\in D_i(m_{\sigma_j^{(r)}})$, we have $u_{i,r+2lr_i}(m')\neq 0\Rightarrow \exists j, l\in \overline{\sigma_j^{(r)}}$. In particular the given of $m'$ suffices to determine the $\overline{\sigma_j^{(r)}}$ : for example we can find the set $\mathcal{M}=\{\text{max}(\overline{\sigma_j^{(r)}})/j,r\}$ and $\mathcal{M}'=\{\text{min}(\overline{\sigma_j^{(r)}})/j,r\}$ in the following way :

if $u_{i,l}(m')=1$ and $u_{i,l+2r_i}(m')=0$ and $u_{i,l+4r_i}(m')\geq 0$, then $l+2r_i\in\mathcal{M}$

if $u_{i,l}(m')=-1$ and $u_{i,l+2r_i}(m')\geq 0$, then $l\in\mathcal{M}$

if $u_{i,l}(m')=-1$ and $u_{i,l-2r_i}(m')=0$ and $u_{i,l-4r_i}(m')\leq 0$, then $l-2r_i\in\mathcal{M}'$

if $u_{i,l}(m')=1$ and $u_{i,l-2r_i}(m')\leq 0$, then $l\in\mathcal{M}'$

\noindent So if $m'\in D_i(m_1)\cap D_i(m_2)$, we have the same decomposition for $m_1$ and $m_2$, that is to say $m_1=m_2$.\qed

\begin{lem}\label{descent} Let $m\in B$, $m'\in D(m)\cap B_j$. We suppose that for all $m''\in D(m)$ such that $v(m''m^{-1})< v(m'm^{-1})$, all $i\in I, l\in\ZZ$ we have $u_{i,l}(m'')\leq 1$. Then $D_j(m')\subset D(m)$.\end{lem}

\demo Let $p=v(m'm^{-1})$. Consider the decomposition of proposition \ref{jdecomp} with $J=\{j\}$ : $L(m)=\underset{M\in B_j\cap D(m)}{\sum}\lambda_j(M)L_j(M)$. It follows from the hypothesis and from the lemma \ref{domiseul} that for $v(Mm^{-1})<p$, $m'\notin D_j(M)$. So $\lambda_j(m')>0$, and $D_j(m')\subset D(m)$.\qed

\begin{prop}\label{cpfacile} Let $i\in I$ such that all $m\in D(Y_{i,L})$ satisfies : for $j\in I$, if $m\in B_j$ then $\forall l\in\ZZ, u_{j,l}(m)\leq 1$. Then all coefficients of $L(Y_{i,L})$ are equal to 1.\end{prop}

\demo We can compute the coefficients of $L(Y_{i,L})=F(Y_{i,L})$ thanks to the classical algorithm (see section \ref{notadeux}) : let us show by induction on $v(mY_{i,L}^{-1})$ that the coefficients of $m$ is equal to $1$. For a monomial $m < Y_{i,L}$, there is $j\in I$ such that $m\notin B_j$. There is $M\rightarrow_j m$. It follows from the lemma \ref{inter} that $M$ is entirely determined by $m$. So the coefficient of $m$ is the coefficient of $M$ in $L(Y_{i,l})$ multiplied with the coefficient of $m$ in $L_j(M)=F_j(M)$, that is to say 1 (section \ref{notadeux}).\qed

\section{Type $A$, $A^{(1)}$}\label{typea} 

\begin{prop}\label{cpaconj} The property of theorem \ref{mainprem} is true for $\Glie$ of type $A_n$ ($n\geq 1$) or $A_l^{(1)}$ ($l\geq 2$).\end{prop}

\noindent Technical consequences of this result which will be used in the following are discussed in section \ref{consa}.

\subsection{Type $A$}

Let $n\geq 1$ and $\Glie$ of type $A_n$. For $i\in\{2,...,n-1\}$, $l\in\ZZ$ :
$$A_{i,l}=Y_{i,l+1}Y_{i,l-1}Y_{i+1,l}^{-1}Y_{i-1,l}^{-1}$$
$$A_{1,l}=Y_{1,l+1}Y_{1,l-1}Y_{2,l}^{-1}\text{ , }A_{n,l}=Y_{n,l+1}Y_{n,l-1}Y_{n-1,l}^{-1}$$
In particular for all $i\in I, l\in\ZZ$, $u(A_{i,l}^{-1})\leq 0$. So $m\leq m'\Rightarrow u(m)\leq u(m')$. 

\noindent We can suppose $Y_{i,L}=Y_{i,0}$ (proposition \ref{aidafm}).

\begin{lem}\label{argu} For $m\in B$ and $m'\in D(m)$ we have $u(m)\geq v_n(m'm^{-1})$.\end{lem}

\demo For all $i\in I$, we have $\omega_{i}+\omega_{n+1-i}\in \alpha_n + \underset{j\leq n-1}{\sum}\ZZ \alpha_j$ (see \cite{bou}). 

\noindent Consider $m'\in D(m)$. It follows from the lemma 6.8 of \cite{Fre2} that $\omega(m)\geq \omega(m')\geq -\underset{i\in I}{\sum}u_i(m)\omega_{n+1-i}$, and so $-\omega(m'm^{-1})\leq \underset{i\in I}{\sum}u_i(m)(\omega_i+\omega_{n+1-i})$. So $v_n(m'm^{-1})\leq \underset{i\in I}{\sum}u_i(m)=u(m)$.\qed

\begin{lem}\label{aidean} For $j\in I$, if $m\in B_j\cap D(Y_{i,0})$ then $u_j(m)\leq 1$.\end{lem}

\demo Suppose there is $j\in I$ and $m_1\in B_j\cap D(Y_{i,0})$ such that $u_j(m_1)\geq 2$. Let $J_1=\{k\in I/k<j\}$, $J_2=\{k\in I/k > j\}$ and $J=J_1\cup J_2$. Let $m_2\rightarrow_J m_1$ and $v=v_{j-1}(m_1m_2^{-1})+ v_{j+1}(m_1m_2^{-1})$. It follows from lemma \ref{argu} (for $\Glie_{J_1}$ and $\Glie_{J_2}$) that $u_{J_1}(m_2)+u_{J_2}(m_2)\geq v$. Moreover we have $u_j(m_2)=u_j(m_1)-v\geq 2-v$. So $u(m_2)=u_{J_1}(m_2)+u_j(m_2)+u_{J_2}(m_2)\geq 2$, contradiction because $m_2\leq Y_{i,0}$.\qed

\noindent The proposition \ref{cpaconj} for type $A_n$ follows from proposition \ref{cpfacile} and lemma \ref{aidean}.

\subsection{Type $A_l^{(1)}$} Let $l\geq 2$ and $\Glie$ of type $A_l^{(1)}$. For $i\in I$, $l\in\ZZ$ (where $Y_{-1,L}=Y_{n,L}$, $Y_{n+1,L}=Y_{0,L}$):
$$A_{i,l}=Y_{i,l+1}Y_{i,l-1}Y_{i+1,l}^{-1}Y_{i-1,l}^{-1}$$
In particular for all $i\in I, l\in\ZZ$, $u(A_{i,l}^{-1})\leq 0$. So $m\leq m'\Rightarrow u(m)\leq u(m')$. 

\noindent We have an analog of lemma \ref{aidean} by putting in the proof $J=I-\{j\}$ instead of $J_1\cup J_2$. In particular we get proposition \ref{cpaconj} for type $A_l^{(1)}$.

\subsection{Consequences}\label{consa} In this section $\Glie$ is general and consider $J\subset I$ such that $\Glie_J$ is of type $A_m$ ($m\leq n$). We prove technical results which will be useful in the following. Let $i\in I, a\in\CC^*$. 

\begin{lem}\label{argudeux} Let $m\in B_J$, $j\in J$ and $m'\in B_j$ such that $m'\in D_J(m)$. We have $u_J(m)\geq u_j(m')$.\end{lem}

\demo It follows from lemma \ref{prodstan} that we can write :
$$m'=m^{(I-J)}\underset{k\in J,l\in\ZZ}{\prod} m'_{k,l,1}...m'_{k,l,u_{k,l}(m)}$$ 
where $m'_{k,l,1}\in D_J(Y_{k,l})$. For $\alpha\in J\times\ZZ\times\NN$, it follows from lemma \ref{aidean} that $u_{j,L}^+(m'_{\alpha})\geq 1\Rightarrow (m'_{\alpha})^{(j)}=Y_{j,L}$. So there are $\alpha_1,...,\alpha_{u_j(m')}$ such that the $(m'_{\alpha_p})^{(j)}=Y_{j,l_p}$ and $m'_{\alpha_1}...m'_{\alpha_{u_j(m')}}=(m')^{(j)}$. So $u_j(m')\leq \underset{k\in J,l\in\ZZ}{\sum}u_{k,l}(m)=u_J(m)$.\qed

\begin{lem}\label{argutrois} Let $M\in B_J$ such that $u_J(M)\geq 2$. The following properties are equivalent :

(i) there are $j\in J$, $l\in\ZZ$, $M_1\in D_J(M)\cap B_j$ such that $u_{j,l}(M_1)\geq 2$

(ii) there are $M'\in D_J(M)\cap B_J$, $i_1,i_2\in J$, $l_1,l_2\in\ZZ$ such that $i_2-i_1\geq |l_1-l_2|$ and $(i_2-i_1)-(l_2-l_1)$ is even and $u_{i_1,l_1}(M')\geq 1$, $u_{i_2,l_2}(M')\geq 1$.

\noindent Moreover one can choose $M'$ such that $M_1\in D_J(M')$.\end{lem}

\demo We can suppose that $\Glie=\Glie_J$ is of type $A_n$. For $K\subset I$, in this proof the notation $\rightarrow_K, \rightharpoonup_K$ is defined as in definition \ref{fleche} by putting $D(M)$ instead of $D(Y_{i,a})$.

\noindent Let us show that $(ii)\Rightarrow (i)$ : if $i_2=i_1$ we have $u_{i_1, l_1}(M')\geq 2$. If $i_2-i_1>0$, suppose that $(i)$ is not true. In this situation we can use the lemma \ref{descent}. Consider the integers:
$$K=\frac{(i_2-i_1)+(l_2-l_1)}{2}\text{ , }K'=\frac{(i_2-i_1)+(l_1-l_2)}{2}$$
We have $K,K'\geq 0$. Denote $i=i_1+K=i_2-K'$, $l=l_1+K=l_2+K'$ and consider :
$$V=A_{i_1,l_1+1}^{-1}A_{i_1+1, l_1+2}^{-1}...A_{i_1+(K-1), l_1+K}^{-1}A_{i_2,l_2+1}^{-1}A_{i_2-1,l_2+2}^{-1}...A_{i_2-(K'-1),l_2+K'}^{-1}$$ 
There is $M_1\in D(M')$ such that $M_1\leq M'V$ and $v_i(M_1(M')^{-1})=0$ (lemma \ref{descent}). In particular $M_1\in B_i$ and $u_{i,l}(M_1)\geq 2$, contradiction.

\noindent Let us show that $(i)\Rightarrow (ii)$ : it follows from lemma \ref{prodstan} and proposition \ref{cpaconj} that we can suppose that $u(M)=2$. Denote $M=Y_{i_1,l_1}Y_{i_2,l_2}$, and let us show the result by induction on $n$. For $n=1$ we have $M_1=M$ and $(ii)$ is clear. In general let $M_1\in D(M)\cap B_j$ such that $M_1^{(j)}=Y_{j,l}^2$. We can suppose that $v(M_1M^{-1})$ is minimal. If $M_1$ is dominant, we put $M_1=M'$. Otherwise consider $J'=\{1,...,n-1\}$ if $j\leq n-1$, and ${J'}=\{2,...,n\}$ if $j=n$ (we suppose that $j\leq n-1$, the case $j=n$ can be treated in the same way). Let $M_2\rightarrow _{J'} M_1$. The induction with $\Glie_{J'}$ of type $A_{n-1}$ gives that $M_2^{({J'})}=Y_{i_1,l_1}Y_{i_2,l_2}$ where $i_2-i_1\geq |l_1-l_2|$ and $(i_2-i_1)-(l_2-l_1)$ is even. We have $u(M_2)\leq 2$ and so $u_n(M_2)=u(M_2)-u_{J'}(M_2)\leq 0$. If $M_2^{(n)}=1$, we put $M_2=M'$. Otherwise it follows from the lemma \ref{argudeux} that we are in one the following cases $\alpha, \beta, \gamma$:

$\alpha$) if $M_2^{(n)}=Y_{n,K_1}^{-1}Y_{n,K_2}^{-1}$, we have : 
$$M_2\leftarrow_n M_3=Y_{i_1,l_1}Y_{i_2,l_2}Y_{n-1,K_1-1}^{-1}Y_{n-1,K_2-1}^{-1}Y_{n,K_1-2}Y_{n,K_2-2}$$ 
If $Y_{i_2,l_2}\neq Y_{n-1,K_1-1}$ and $Y_{i_2,l_2}\neq Y_{n-1,K_2-1}$, there is $M_4=M_3(M_1M_2^{-1})\in D(M_3)$ (lemma \ref{descent}) such that $M_4^{(j)}=Y_{j,l}^2$ and $v(M_4M^{-1})<v(M_1M^{-1})$, contradiction. So for example we have $M_3=Y_{i_1,l_1}Y_{n-1,K_2-1}^{-1}Y_{n,K_1-2}Y_{n,K_2-2}$ and $i_2=n-1$, $l_2=K_1-1$. We have : 
$$M_3\leftarrow_{\{i_1+1,...,n-1\}} M_4=Y_{i_1,l_1}Y_{i_1,K_2-1-(n-i_1-1)}^{-1}Y_{i_1+1,K_2-2-(n-i_1-1)}Y_{n,K_1-2}$$ 
If $Y_{i_1,l_1}Y_{i_1,K_2-1-(n-i_1-1)}^{-1}\neq 1$, there is :
$$M_4\leftarrow_{\{1,...,i_1\}} M_5=Y_{1,K_3}Y_{i_1,l_1}Y_{i_1+1,K_2-2-(n-i_1-1)}^{-1}Y_{i_1+1,K_2-2-(n-i_1-1)}Y_{n,K_1-2}$$ 
and $u(M_5)=3$, impossible. So $l_1=K_2-1-(n-i_1-1)$ and $M_4=Y_{i_1',l_1'}Y_{i_2',l_2'}$ where $i_1'=i_1+1$, $i_2'=n$, $l_1'=K_2-2-(n-i_1-1)$, $l_2'=K_1-2$. Let $M'=M_4$. We have $i_2'-i_1'=i_2-i_1\geq |l_2-l_1|=|l_2'-l_1'|$ and $(i_2'-i_1')-(l_2'-l_1')=(i_2-i_1)-(l_2-l_1)$ is even.

$\beta$) if $M_2^{(n)}=Y_{n,K_1}Y_{n,K_2}^{-1}$, we have :
$$M_2\leftarrow_n M_3=Y_{i_1,l_1}Y_{i_2,l_2}Y_{n-1,K_2-1}^{-1}Y_{n,K_1-2}Y_{n,K_2-2}$$
and $u(M_3)=3$, impossible.

$\gamma$) if $M_2^{(n)}=Y_{n,K_1}^{-1}$, we have : 
$$M_2\leftarrow_n M_3=Y_{i_1,l_1}Y_{i_2,l_2}Y_{n-1,K_1-1}^{-1}Y_{n,K_1-2}$$ 
If $Y_{n-1,K_1-1}\neq Y_{i_1,l_1}$ and $Y_{n-1,K_1-1}\neq Y_{i_2,l_2}$, there is $M_4=M_3(M_1M_2^{-1})\in D(M_3)$ (lemma \ref{descent}) such that $M_4^{(j)}=Y_{j,l}^2$ and $v(M_4M^{-1})<v(M_1M^{-1})$, contradiction. So for example we have $M_3=Y_{i_1,l_1}Y_{n,K_1-2}$ and $i_2=n-1$, $l_2=K_1-1$. Let $M'=M_3$ and we have $n-i_1=i_2-i_1+1\geq |l_2-l_1|+1\geq |(K_1-2)-l_1+1|+1\geq |(K_1-2)-l_1|$ and $n-i_1-((K_1-2)-l_1)=(i_2-i_1)-(l_2-l_1)$ is even.

\noindent For the last point, the arguments of this proof can be used for any $M'\in B$ such that $M_1\in D_J(M')$.\qed

\section{Type $B$}\label{typec}

\subsection{Statement} In this section $\Glie$ is of type $B_n$ ($n\geq 2$). For $i\in\{2,...,n-2\}$, $l\in\ZZ$ :
$$A_{i,l}=Y_{i,l+2}Y_{i,l-2}Y_{i+1,l}^{-1}Y_{i-1,l}^{-1}\text{ , }A_{1,l}=Y_{1,l+2}Y_{1,l-2}Y_{2,l}^{-1}\text{ , }A_{n,l}=Y_{n,l+1}Y_{n,l-1}Y_{n-1,l}^{-1}$$
$$A_{n-1,l}=Y_{n-1,l+2}Y_{n-1,l-2}Y_{n-2,l}^{-1}Y_{n,l-1}^{-1}Y_{n,l+1}^{-1}$$

\noindent In this section we prove: 

\begin{prop}\label{resbb} The property of theorem \ref{mainprem} is true for $\Glie$ of type $B_n$ ($n\geq 2$).\end{prop}

\noindent Denote  $J=\{1,...,n-1\}$. We can suppose $Y_{i,L}=Y_{i,0}$ (proposition \ref{aidafm}). 

\noindent As $u(A_{n-1,l}^{-1})>0$, the $m'\leq m$ does not imply $u(m')\leq u(m)$. 

\subsection{Proof of the proposition \ref{resbb}}

Suppose that there is $m\in D(Y_{i,0})$ such that there is $j\in J, l\in\ZZ$, $u_{j,l}(m)\geq 2$, and let $m$ such that $v(mY_{i,0}^{-1})$ is minimal with this property.

\begin{lem}\label{ndessus}  There is $M\in D(Y_{i,0})$ such that $v(MY_{i,0}^{-1})<v(mY_{i,0}^{-1})$ and $\exists l'\in\ZZ$, $u_{n,l'}(M)\geq 2$.\end{lem}

\demo Suppose that $M$ does not exists. Let $m_1\rightarrow_J m$. It follows from lemma \ref{argutrois} that $m_1=m_1'Y_{i_1,l_1}Y_{i_2,l_2}$ where $m_1'\in B_J$, $2(i_2-i_1)\geq |l_1-l_2|$, $(i_2-i_1)-(l_2-l_1)/2$ is even. Let $m_2\rightharpoonup_n m_1$ such that $v(m_1m_2^{-1})=1$ ($m_2$ exists because it follows from the hypothesis and lemma \ref{descent} that $M_2\rightarrow_n m_2$ implies $D_n(M_2)\subset D(Y_{i,0})$). We have $m_2=m_2'Y_{i_1,l_1}Y_{i_2,l_2}Y_{n-1,L}^{-1}Y_{n,L-1}$ where $m_2'\in B_J$ and $u_{n,L-1}(m_2')\geq 0$. If $Y_{i_1,l_1}Y_{i_2,l_2}Y_{n-1,L}^{-1}\notin B$, there is $m_2\rightharpoonup_J M_2=m_2 (mm_1^{-1})$ (lemma \ref{descent}) such that $u_{j,l}(M_2)\geq 2$ and $v(M_2Y_{i,0}^{-1})<v(mY_{i,0}^{-1})$, contradiction. So for example $Y_{i_2,l_2}=Y_{n-1,L}$, and $m_2\in B_J$. $m_2$ is not dominant (because we would have $u(m_2)\geq 2$), so there is $m_3\rightharpoonup_n m_2$ such that $v(m_2m_3^{-1})=1$ (same argument as above for the existence of $m_3$). We have $m_3=m_3'Y_{i_1,l_1}Y_{n-1,L'}^{-1}Y_{n,L-1}Y_{n,L'-1}$ where $m_3'\in B_J$ and $u_{n,L-1}(m_3')\geq 0$, $u_{n,L'-1}(m_3')\geq 0$. 

if we can choose $L'\neq L+2$, the same argument gives $Y_{i_1,l_1}=Y_{n-1,L'}$ because :
$$m_4=Y_{i_1,l_1}Y_{i_2,l_2}Y_{n-1,L'}^{-1}Y_{n,L-1'}m_4'\rightharpoonup_n m_1$$ where $m_4'\in B_J$. So we have $i_1=i_2=n-1$, so $l_1=l_2$ and $L'=L$, ie $u_{n,L}(m_3)\geq 2$, contradiction.

if we can not choose $L'\neq L+2$, we can not use the same argument (because : $1\notin L_n(Y_{n,L-1}Y_{n,L+1})$). We have ($k\geq 1$):
$$m_3\leftarrow_n m_5=m_5'Y_{i_1,l_1}Y_{n-1,L+2}^{-1}Y_{n-1,L+4}^{-1}...Y_{n-1,L-1+2k+1}^{-1}Y_{n,L-1}Y_{n,L+1}...Y_{n,L-1+2k}$$
where $m_5'\in B$. Suppose that $m_5'Y_{i_1,l_1}Y_{n-1,L+4}^{-1}\notin B$. Then we have : 
$$m_5\leftharpoonup_{n-1} m_6=m_5'Y_{i_1,l_1}Y_{n-1,L}Y_{n-1,L+2}^{-1}Y_{n-1,L+6}^{-1}...Y_{n-1,L-1+2k+1}^{-1}Y_{n,L-1}Y_{n,L+5}...Y_{n,L-1+2k}$$ 
But $u_{n-1,L}(m_6A_{n,L}^{-1})=2$, contradiction. In the same way we prove by induction that 
\\$m_5'Y_{i_1,l_1}Y_{n-1,L+4}^{-1}Y_{n-1,L+8}^{-1}...\in B$ and so : 
$$m_5=m_5''Y_{i_1,l_1}Y_{n-1,L+2}^{-1}Y_{n-1,L+6}^{-1}...Y_{n-1,L+2+4K'}^{-1}Y_{n,L-1}Y_{n,L+1}...Y_{n,L-1+2k}$$
where $m_5''\in B$. Suppose that $m_5\notin B$. As $m_5$ is right negative, we have $L+2+4K'=L+2k$ and so $k=1+2K'$. Consider $m_7\rightarrow_J m_5$. We get that $m_7$ is dominant, and so $m_7=Y_{i,0}$. Let $K''=u_-(m_5)\leq K'$. We have $1=u(m_7)\geq u_J(m_7)+u_n(m_7)\geq K''+(k+1-2K'')=1+k-K''\geq k-K'\geq 1+K'$. So $K'=0$ and $k=1$. In particular $m_5=m_5''Y_{i_1,l_1}Y_{n-1,L+2}^{-1}Y_{n,L-1}Y_{n,L+1}$ and so we have $i_1=n-1-j$, $l_1=L+2-2j$ where $j\geq 0$ (otherwise we would have $u(m_7)\geq 2$). It implies $i_2-i_1-(l_2-l_1)/2=j-(j-1)=1$ not even, contradiction. So $m_5\in B$ and $m_5=Y_{i,0}$. But $u(m_5)\geq u_n(m_5)\geq 2$, contradiction. \qed

\begin{lem}\label{capcdeux} Let $j\in I$ and $m\in D(Y_{i,0})\cap B_j$ such that $u_j(m)=2$. For $L,L'\in\ZZ$ such that $m^{(j)}=Y_{j,L}Y_{j,L'}$, we have $L\neq L'$.\end{lem}

\demo It follows from lemma \ref{ndessus} that we can suppose that $j=n$ and that for $v(m'Y_{i,0}^{-1})\leq v(mY_{i,0}^{-1})$, for all $j\in J, l\in\ZZ$, $u_{j,l}(m')\leq 1$. Let $M$ such that $\exists l\in\ZZ$, $u_{n,l}(M)\geq 2$ and suppose that $v(MY_{i,0}^{-1})$ is minimal with this property. Let $L$ be maximal such that $u_{n,L-1}(M)\geq 2$. First it follows from lemma \ref{argutrois} that $M\in B_n$. We have $M\rightarrow_n M'=MA_{n,L}^{-1}$ and the coefficient of $M'$ in $L_n(M)$ is at least 2. Suppose that there is $j\in J$ such that $M'\notin B_j$. Let $M''\rightarrow_j M'$. It follows from lemma \ref{inter} that $M''$ is uniquely determined by $M'$, and that the coefficient of $M'$ in $F_j(M')$ is $1$. But the coefficient of $M''$ in $L(Y_{i,0})$ is 1, so it follows from the proposition \ref{jdecomp} that the coefficient of $M'$ is 1, contradiction. So $M'\in B_J$. So $M=Y_{n-1,L}^{-1}\tilde{M}$ where $\tilde{M}\in B$ and $u_{n,L-1}(\tilde{M})\geq 2$. As $u_n(M)\geq 2$, $M\notin B$. So there is $M_0\rightarrow_J M$. We have :
$$M_0=\tilde{M_0}Y_{n,L-1}Y_{n,L-3}^{-1}$$
where $\tilde{M_0}\in B$ and $u_J(\tilde{M_0})\geq 1$. But $M_0$ is not right negative, so $M_0$ is dominant (corollary \ref{propgene}). But $u(M_0)\geq u_J(\tilde{M}_0)+u_{n,L-1}(M_0)\geq 2$, contradiction. \qed

\noindent So the proposition \ref{resbb} follows from proposition \ref{cpfacile} and lemma \ref{capcdeux}.

\subsection{Complement : degree of monomials} 

The aim of this section is to prove that the degrees are bounded (it is a complement independent of the proof of theorem \ref{mainprem}):

\begin{prop}\label{born} For $j\in I$ and $m\in B_j\cap D(Y_{i,0})$, then $u_j(m)\leq 2$.\end{prop}

\noindent Note that it follows from proposition \ref{resbb} that we can use the lemma \ref{descent}.

\noindent For $m\in A$ denote $w(m)=(u_J^+(m),u_J^-(m),u_n^+(m),u_n^-(m))$.

\noindent Suppose that there is $j\in J$ and $m_0\in D(Y_{i,0})\cap B_j$ such that $u_j(m_0)\geq 3$. It follows from lemma \ref{argudeux} that there is $m\rightarrow_J m_0$ such that $u_J(m)\geq 3$. Suppose that $v(mY_{i,0}^{-1})$ is minimal for this property.

\begin{lem}\label{commence} There is $M\in D(Y_{i,0})\cap B_n$ such that $M>m$ and $u_n(M)\geq 3$.\end{lem}

\demo We have $m\in B_J$ and $u_J(m)=3$ : $m=Y_{i_1,l_1}Y_{i_2,l_2}Y_{i_3,l_3}m'$ where $(m')^{(J)}=1$. There is 
$$m\leftharpoonup_n m_1=Y_{i_1,l_1}Y_{i_2,l_2}Y_{i_3,l_3}Y_{n-1,L}^{-1}m_1'Y_{n,L-1}$$ where $(m_1')^{(J)}=1$
and $u_{n,L-1}(m_1')\geq 0$. If $Y_{i_1,l_1}Y_{i_2,l_2}Y_{i_3,l_3}Y_{n-1,L}^{-1}\notin B$, there is $M_1\rightarrow_J m_1$ such that $u_J(M_1)\geq 3$, contradiction. So for example $m_1=Y_{i_1,l_1}Y_{i_2,l_2}m_1'Y_{n,l_3-1}$. There is 
$$m_1\leftharpoonup_n m_2=Y_{i_1,l_1}Y_{i_2,l_2}Y_{n-1,L'}^{-1}m_2'Y_{n,l_3-1}Y_{n,L'-1}$$ 
where $(m_2')^{(J)}=1$ and $u_{n,l_3-1}(m_2')\geq 0$, $u_{n,L'-1}(m_2')\geq 0$. It follows from lemma \ref{capcdeux} that $l_3\neq L'$. If $L'=l_3+1$, we have $m_2''\rightarrow_J m_2$ where $m_2''$ is dominant and $u(m_2'')\geq u_J(m_2'')\geq 2$, contradiction. So we see as above that for example $m_2=Y_{i_1,l_1}m_2'Y_{n,l_3-1}Y_{n,l_2-1}$. In the same way 
$$m_2\leftharpoonup_n m_3=m_3'Y_{i_1,l_1}Y_{n-1,L''}^{-1}Y_{n,l_3-1}Y_{n,l_2-1}Y_{n,L''-1}$$ 
where $(m_3')^{(J)}=1$ and $u_{n,l_3-1}(m_3')\geq 0$, $u_{n,l_2-1}(m_3')\geq 0$, $u_{n,l_1-1}(m_3')\geq 0$. We can conclude with lemma \ref{argudeux}.\qed

\noindent For $m\in A$, denote $w(m)=(u_J^+(m),u_J^-(m),u_n^+(m),u_n^-(m))$.

End of the proof of proposition \ref{born} :

\noindent Suppose that there is $M\in D(Y_{i,0})\cap B_n$ such that $u_n(M)\geq 3$. Suppose that $v(MY_{i,0}^{-1})$ is minimal for this property. It follows from lemma \ref{commence} that for $M'\in D(Y_{i,0})$, $M'\geq M$ for $j\in J, l\in\ZZ$, $u_{j,l}(M')\leq 2$.

\noindent We have $M\in B_{\{1,...,n-2\}}$ (if not we would have $M'\rightarrow_{\{1,...,n-2\}} M$ with $M'>M$ and $u_n(M')\geq 3$).

If $u_n(M)=3$ : there is $M\leftarrow_J M_1$. If $v_{n-1}(M_1M^{-1})=1$, we have $w(M_1)=(a,0,3,2)$ or $(a,0,2,1)$ or $(a,0,1,0)$ where $a=1$ or $2$. For the first two cases we have $u_n^+(M_1)+u_n^-(M_1)\geq 3$, so there is $M_1'\rightarrow_n M_1$ such that $u_n(M_1)\geq 3$, contradiction. For the last case $M_1$ is dominant with $u(M_1)\geq 2$, contradiction. So $v_{n-1}(M_1M^{-1})=2$ and $w(m_1)=(2,0,3,4)$ or $(2,0,2,3)$ or $(2,0,1,2)$ or $(2,0,0,1)$. As above we have $w(M_1)=(2,0,0,1)$. Let $M_2\rightarrow_n M_1$. If $w(M_2)=(1,0,1,0)$, $M_2$ is dominant with $u(M_2)\geq 2$, contradiction. So $w(M_2)=(2,1,1,0)$. Let $M_3\rightarrow_J M_2$. If $w(M_3)=(2,0,1,2)$, there is $M_3'>M_3$ such that $u_n(M_3')\geq 3$, contradiction. So $w(M_3)=(2,0,0,1)$. We continue and we get an infinite sequence such that $w(M_{2k})=(2,1,1,0)$ and $w(M_{2k+1})=(2,0,0,1)$. Contradiction because the sequence $v(M_kY_{i,0}^{-1})\geq 0$ decreases strictly.

If $u_n(M)=4$ : there is $M\leftarrow_J M_1$. If $v_{n-1}(M_1M^{-1})=1$, we have $w(M_1)=(a,0,4,2)$ or $(a,0,3,1)$ or $(a,0,2,0)$ where $a=1$ or $2$. We see as above that it is impossible. So $v_{n-1}(M_1M^{-1})=2$ and $w(m_1)=(2,0,4,4)$ or $(2,0,3,3)$ or $(2,0,2,2)$ or $(2,0,1,1)$. As above we have $w(M_1)=(2,0,1,1)$. Let $M_2\rightarrow_n M_1$. If $w(M_2)=(1,0,2,0)$, $M_2$ is dominant with $u(M_2)\geq 2$, contradiction. So $w(M_2)=(2,1,2,0)$. Let $M_3\rightarrow_J M_2$. If $w(M_3)=(2,0,2,2)$, there is $M_3'>M_3$ such that $u_n(M_3')\geq 3$, contradiction. So $w(M_3)=(2,0,1,1)$. We continue and we get an infinite sequence such that $w(M_{2k})=(2,1,2,0)$ and $w(M_{2k+1})=(2,0,1,1)$. Contradiction because the sequence $v(M_kY_{i,0}^{-1})\geq 0$ decreases strictly.\qed

\section{Type $C$}\label{typeb}

\subsection{Statement} Let $\Glie$ be of type $C_n$ ($n\geq 2$). For $i\in\{2,...,n-1\}$, $l\in\ZZ$ :
$$A_{i,l}=Y_{i,l+1}Y_{i,l-1}Y_{i-1,l}^{-1}Y_{i+1,l}^{-1}$$
$$A_{1,l}=Y_{1,l+1}Y_{1,l-1}Y_{2,l}^{-1}\text{ , }A_{n,l}=Y_{n,l+2}Y_{n,l-2}Y_{n-1,l+1}^{-1}Y_{n-1,l-1}^{-1}$$
In particular for all $i\in I, l\in\ZZ$, $u(A_{i,l}^{-1})\leq 0$. So $m\leq m'\Rightarrow u(m)\leq u(m')$.

\noindent In this section we prove :

\begin{prop}\label{rescc} The property of theorem \ref{mainprem} is true for $\Glie$ of type $C_n$ ($n\geq 2$).\end{prop}

\noindent Denote $J=\{1,...,n-1\}\subset I$. 

\subsection{Proof of proposition \ref{rescc}} We can suppose $Y_{i,L}=Y_{i,0}$ (proposition \ref{aidafm}). 

\begin{lem}\label{capb}  (i) For $m\in B_n\cap D(Y_{i,0})$, we have $u_n(m)\leq 1$.

(ii) For $j\leq n-1$ and $m\in B_j\cap D(Y_{i,0})$, we have $u_j(m)\leq 2$. 

(iii) Let $j\leq n-1$ and $m\in D(Y_{i,0})\cap B_j$ such that $u_j(m)=2$. For $L,L'\in\ZZ$ such that $m^{(j)}=Y_{j,L}Y_{j,L'}$, we have $L\neq L'$.\end{lem}

\demo $(i)$ suppose that there is $m_1\in B_n\cap D(Y_{i,0})$ such that $u_n(m_1)\geq 2$. Let $m_2\rightarrow_J m_1$. We have $u_n(m_2)\geq 2-v_{n-1}(m_2m_1^{-1})$. But $u_J(m_2)\geq v_{n-1}(m_2m_1^{-1})$ (lemma \ref{argu}) and so $u(m_2)=u_J(m_2)+u_n(m_2)\geq 2$. As $Y_{i,0}\geq m_2$ it is impossible.

\noindent $(ii)$ suppose that there is $j\leq n-1$ and $m_1\in B_j\cap D(Y_{i,0})$ such that $u_j(m_1)\geq 3$. Let $m_2\rightarrow_J m_1$. Then we have $u_J(m_2)\geq 3$ (lemma \ref{argudeux}) and so $u(m_2)\geq 3+u_n(m_1)\geq 2$ (it follows from $(i)$ that $u_n(m_1)\geq -1$). Contradiction.

\noindent $(iii)$ let $j\neq n$ and $m_1\in D(Y_{i,0})\cap B_j$ such that $m_1^{(j)}=Y_{j,L}^2$. We can suppose that $v(m_1m^{-1})$ is minimal. Let $m_2\rightarrow_J m_1$. It follows from lemma \ref{argutrois} for $\Glie_J$ of type $A_{n-1}$ that $m_2^{(J)}=Y_{i_1,L_1}Y_{i_2,L_2}$ with $i_2-i_1\geq |L_1-L_2|$ and $(i_2-i_1)-(L_2-L_1)$ is even. As $u(m_2)\leq 1$ and $u_n(m_2)\geq -1$, we have $u_n(m_2)=-1$ and $u_J(m_2)=2$. In particular $m_2=Y_{i_1,L_1}Y_{i_2,L_2}Y_{n,K}^{-1}$. There is $m_2\leftarrow_n m_3=Y_{i_1,L_1}Y_{i_2,L_2}Y_{n-1,K-1}^{-1}Y_{n-1,K-3}^{-1}Y_{n,K-4}$. We are in one the following cases $\alpha, \beta, \gamma, \delta$:

$\alpha$) $Y_{i_1,L_1}Y_{i_2,L_2}Y_{n-1,K-1}^{-1}Y_{n-1,K-3}^{-1}=1$ : impossible because $i_1=i_2\Rightarrow L_1=L_2$. 

$\beta$) $Y_{i_1,L_1}Y_{i_2,L_2}Y_{n-1,K-1}^{-1}Y_{n-1,K-3}^{-1}=Y_{i_1,L_1}Y_{n-1,K-1}^{-1}\neq 1$ (or in the same way $Y_{i_2,L_2}Y_{n-1,K-1}^{-1}\neq 1$). There is $m_3\leftarrow_{n-1} m_4=m_3A_{n-1,K-2}$. In particular $m_4^{(n)}=Y_{n,K-4}Y_{n,K-2}^{-1}$, contradiction with $(i)$.

$\gamma$) $Y_{i_1,L_1}Y_{i_2,L_2}Y_{n-1,K-1}^{-1}Y_{n-1,K-3}^{-1}=Y_{i_1,L_1}Y_{n-1,K-3}^{-1}\neq 1$ (or in the same way $Y_{i_2,L_2}Y_{n-1,K-3}^{-1}\neq 1$). In particular $i_2=n-1$, $L_2=K-1$ and $m_3=Y_{i_1,L_1}Y_{n-1,K-3}^{-1}Y_{n,K-4}$. Let $J'=\{i_1+1,...,n-1\}$ ($J'$ can be empty) and $m_3\leftarrow_{J'} m_4=Y_{i_1,L_1}Y_{i_1,K-2-n+i_1}^{-1}Y_{i_1+1,K-3-n+i_1}$. If $Y_{i_1,L_1}Y_{i_1,K-2-n+i_1}^{-1}\neq 1$, let $m_5\rightarrow_{i_1} m_4$. If $i_1=1$, we have $u(m_5)=2$, impossible. If $i_1\geq 2$, we have $m_5=Y_{i_1-1,K-3-n+i_1}^{-1}Y_{i_1,L_1}Y_{i_1,K-4-n+i_1}$. Let $J''=\{1,...,i_1-1\}$ and $m_5\leftarrow_{J''} m_6=Y_{1,K'}Y_{i_1,L_1}$. We have $u(m_6)=2$, impossible. So $Y_{i_1,L_1}=Y_{i_1,K-2-n+i_1}$, that is to say $L_1=K-2-n+i_1$. So $L_2-L_1=n-i_1+1=i_2-i_1+2>i_2-i_1$, contradiction.

$\delta$) $\{(i_1,L_1),(i_2,L_2)\}\cap\{(n-1,K-1),(n-1,K-3)\}$ is empty : there is $m_3\rightarrow_J m_4=m_3(m_1m_2^{-1})$ such that $m_4^{(j)}=Y_{j,l}^2$ (lemma \ref{descent}) and $v(m_4m^{-1})<v(m_1m^{-1})$, contradiction.\qed

\noindent The proposition \ref{rescc} follows from proposition \ref{cpfacile} and lemma \ref{capb}.

\section{Application to $q,t$-characters}\label{qtapp}

In this section we state and prove the main result of this paper on $q,t$-characters (theorem \ref{posqt}).

\subsection{Reminder on $q,t$-characters \cite{Naa, Nab, her01, her02, her03}}\label{remqt}

We define the product $*_t$ on $A\times (\mathcal{A}\otimes \ZZ[t{\pm}])$ such that : for $(m,v), (m',v')\in A\times \mathcal{A}$ ($m,m',v,v'$ monomials):
$$(m,v)*_t (m',v')=t^{D((m,v),(m',v'))}(mm',vv')$$
where :
$$D((m,v),(m',v'))=\underset{i\in I,l\in\ZZ}{\sum}2u_{i,l+r_i}(m)v_{i,l}(v')+2v_{i,l+r_i}(v)u_{i,l}(m')+v_{i,l+r_i}(v)u_{i,l}(v')+u_{i,l+r_i}(v)v_{i,l}(v')$$
(see \cite{Nab} for the $ADE$-case and \cite{her02, her03} for other cases).

\noindent Let $\Yim_t=\Yim\otimes_{\ZZ}\ZZ[t^{\pm}]$. One can define $\mathfrak{K}_{i,t}, \mathfrak{K}_t\subset \Yim_t$ with deformed screening operators (see \cite{her01, her03}).

\begin{defi} We say that a $\ZZ$-linear map $\chi_{q,t}:\text{Rep}\rightarrow\Yim_t$ is a morphism of $q,t$-characters if :

1) For $M$ a $\text{Rep}$-monomial define $m=\underset{i\in I,l\in\ZZ}{\prod}(Y_{i,l})^{x_{i,l}(M)}\in B$. We have :
$$\chi_{q,t}(M)=m+\underset{m'<m}{\sum}a_{m'}(t)m'\text{ (where $a_{m'}(t)\in\ZZ[t^{\pm}]$)}$$

2) The image of $\chi_{q,t}$ is contained in $\mathfrak{K}_t$.

3) Let $M_1,M_2$ be $\text{Rep}$-monomials. If $\text{max}\{l/\underset{i\in I}{\sum}x_{i,l}(M_1)>0\}\leq \text{min}\{l/\underset{i\in I}{\sum}x_{i,l}(M_2)>0\}$ then :
$$(M_1M_2,(M_1M_2)^{-1}\chi_{q,t}(M_1M_2))=(M_1,M_1^{-1}\chi_{q,t}(M_1))*_t(M_2,M_2^{-1}\chi_{q,t}(M_2))$$\end{defi}

\noindent Those properties are generalizations of the axioms that Nakajima \cite{Nab} defined for the $ADE$-case. 

\begin{thm}\label{axiomes}(\cite{Nab, her02, her03}) For $C$ such that $i\neq j\Rightarrow C_{i,j}C_{j,i}\leq 3$, there is a unique morphism of $q,t$-characters.\end{thm}

\noindent This result (among others) was proved by Nakajima \cite{Nab} for $C$ of type ADE. For $C$ finite it is proved in \cite{her02}, and for $C$ such that $i\neq j\Rightarrow C_{i,j}C_{j,i}\leq 3$ in \cite{her03} (it includes quantum affine and toroidal algebras except $A_1^{(1)}, A_2^{(2)}$). The existence of $\chi_{q,t}$ for symmetric toroidal type is also mentioned in \cite{Nad}.

\noindent In \cite{her02} we defined a $t$-deformed algorithm : for $m\in B$, if it is well-defined it gives an element $F_t(m)\in\mathfrak{K}_t$ such that $m$ is the unique dominant monomial of $F_t(m)$ (an algorithm was also used by Nakajima in the $ADE$-case in \cite{Naa}). If we set $t=1$ we get the classical algorithm. It follows from theorem \ref{axiomes} that the $t$-deformed algorithm is well defined if $i\neq j\Rightarrow C_{i,j}C_{j,i}\leq 3$. We proved in \cite{her02} that if the $t$-deformed algorithm is well-defined, for $i\in I,j\in I,l\in\ZZ$ : $F_t(Y_{i,l})F_t(Y_{j,l})=F_t(Y_{j,l})F_t(Y_{i,l})$. 

\noindent Note that $\chi_{q,t}$ is injective and we have (see \cite{her02}):
\begin{equation}\label{defieq}\chi_{q,t}(\underset{i\in I,l\in\ZZ}{\prod}X_{i,l}^{x_{i,l}})=\overset{\rightarrow}{\underset{l\in\ZZ}{\prod}}\underset{i\in I}{\prod}F_t(Y_{i,l})^{x_{i,l}}=...(\underset{i\in I}{\prod}F_t(Y_{i,l-1})^{x_{i,l-1}})(\underset{i\in I}{\prod}F_t(Y_{i,l})^{x_{i,l}})(\underset{i\in I}{\prod}F_t(Y_{i,l+1})^{x_{i,l+1}})...\end{equation}

\subsection{Technical complement}

\begin{prop}\label{fondqt} (i) Let $m\in B_j$ such that for all $l\in\ZZ$, $u_{j,l}(m)\leq 1$. Then $F_{i,t}(m)=F_i(m)=L_i(m)$ and all coefficients are equal to 1.

(ii) Let $i\in I$ such that all $m\in D(Y_{i,L})$ satisfies : for $j\in I$, if $m\in B_j$ then $\forall l\in\ZZ, u_{j,l}(m)\leq 1$. Then $F_t(Y_{i,L})=F(Y_{i,L})=L(Y_{i,L})\in\Yim_t$ is in $\mathfrak{K}_t$ and all coefficients are equal to 1.\end{prop}

\demo (i) Direct consequence of the lemma 4.13 of \cite{her02}.

(ii) Let $j$ be in $I$ and consider the decomposition of proposition \ref{jdecomp} :
$$L(Y_{i,L})=\underset{m'\in B_j\cap D(Y_{i,L})}{\sum}\lambda_j(m') L_j(m')$$
But it follows from (i) that $m'\in B_j\cap D(Y_{i,L})$ implies that $L_j(m')=F_{j,t}(m')$. And so:
$$L(Y_{i,l})=\underset{m'\in B_j\cap D(Y_{i,L})}{\sum}\lambda_j(m')F_{j,t}(m')\in \mathfrak{K}_{j,t}$$
So $L(Y_{i,L})\in\mathfrak{K}_t$ and $F_t(Y_{i,L})=L(Y_{i,L})=F(Y_{i,l})$.\qed

\subsection{New results for $q,t$-characters}\label{newqtres} It follows also from theorem \ref{mainprem} and proposition \ref{fondqt} :

\begin{prop}\label{abcqt} Let $\Glie$ be of type $A_n$ ($n\geq 1$), $A_l^{(1)}$ ($l\geq 2$), $B_{n}$ ($n\geq 2$) or  $C_n$ ($n\geq 2$). For $i\in I, a\in\CC^*$, we have $\chi_{q,t}(V_i(a))=\chi_q(V_i(a))$ and all coefficients are equal to $1$.\end{prop}

\noindent We prove a conjecture of \cite{her02}:

\begin{thm}\label{posqt} Let $\U_q(\hat{\Glie})$ be a quantum affine algebra ($C$ finite) and $M$ be a standard module of $\U_q(\hat{\Glie})$. The coefficients of $\chi_{q,t}(M)$ are in $\NN[t^{\pm}]$ and the monomials of $\chi_{q,t}(M)$ are the monomials of $\chi_q(M)$. \end{thm}

\noindent In particular the $q,t$-characters for quantum affine algebras have a finite number of monomials and this result shows that the $q,t$-characters of \cite{her02} can be considered as a $t$-deformation of $q$-characters for all quantum affine algebras. In particular it is an argument for the existence of a geometric model behind the $q,t$-characters in non simply-laced cases.

\demo If follows from formula (\ref{defieq}) in section \ref{remqt} that it suffices to look at the $F_t(Y_{i,l})$. We do it with a case by case investigation :

the case $ADE$ follows from the work of Nakajima \cite{Nab}

the case $BC$ follows from theorem \ref{mainprem} and proposition \ref{abcqt} (ii)

the case $G_2$ follows from an explicit computation in \cite{her02}

the case $F_4$ follows from an explicit computation on computer (see section \ref{anncalc}).\qed

\section{Appendix : explicit computations on computer for type $F_4$}\label{anncalc}

The proof of theorem \ref{posqt} for type $F_4$ is based on an explicit computation on computer. A computer program written in C with Travis Schedler computes explicitly the $q,t$-characters of fundamental representations.

\noindent For type $F_4$ there are $4$ fundamental representations (see \cite{bou} for the numbers on the Dynkin diagram) : $\text{dim}(V_1(a))=26$ ($26$ monomials), $\text{dim}(V_1(a))=299$ ($283$ monomials), $\text{dim}(V_1(a))=1703$ ($1532$ monomials), $\text{dim}(V_1(a))=53$ ($53$ monomials). We checked that the coefficients are in $\NN[t^{\pm}]$. We give an explicit list of terms of fundamental representations of type $F_4$ whose coefficient is not $1$ (the complete list of monomials can be found on http://www.dma.ens.fr/$\sim$dhernand/f4monomials.pdf). We can see that the coefficient are all $(t+t^{-1})\in\NN[t^{\pm}]$. They appear only in fundamental representations $2$ and $3$ :

Fundamental representation 2 :

\noindent Monomial 70: $ (t^{-1} +t)$ $
Y_{1,10}Y_{2,7}Y^{-1}_{2,9}Y^{-1}_{2,11}Y_{4,6}
$\\Monomial 87: $ (t^{-1} +t)$ $
Y^{-1}_{1,12}Y_{2,7}Y^{-1}_{2,9}Y_{4,6}
$\\Monomial 89: $ (t^{-1} +t)$ $
Y_{1,10}Y_{2,7}Y^{-1}_{2,9}Y^{-1}_{2,11}Y_{3,8}Y^{-1}_{4,10}
$\\Monomial 105: $ (t^{-1} +t)$ $
Y^{-1}_{1,12}Y_{2,7}Y^{-1}_{2,9}Y_{3,8}Y^{-1}_{4,10}
$\\Monomial 109: $ (t^{-1} +t)$ $
Y_{1,10}Y_{2,7}Y^{-1}_{3,12}
$\\Monomial 120: $ (t^{-1} +t)$ $
Y_{1,8}Y_{1,10}Y^{-1}_{2,9}Y_{3,8}Y^{-1}_{3,12}
$\\Monomial 124: $ (t^{-1} +t)$ $
Y^{-1}_{1,12}Y_{2,7}Y_{2,11}Y^{-1}_{3,12}
$\\Monomial 142: $ (t^{-1} +t)$ $
Y_{2,7}Y^{-1}_{2,13}
$\\Monomial 143: $ (t^{-1} +t)$ $
Y_{1,8}Y^{-1}_{1,12}Y^{-1}_{2,9}Y_{2,11}Y_{3,8}Y^{-1}_{3,12}
$\\Monomial 151: $ (t^{-1} +t)$ $
Y^{-1}_{1,10}Y^{-1}_{1,12}Y_{2,11}Y_{3,8}Y^{-1}_{3,12}
$\\Monomial 155: $ (t^{-1} +t)$ $
Y_{1,8}Y^{-1}_{2,9}Y^{-1}_{2,13}Y_{3,8}
$\\Monomial 168: $ (t^{-1} +t)$ $
Y^{-1}_{1,10}Y^{-1}_{2,13}Y_{3,8}
$\\Monomial 173: $ (t^{-1} +t)$ $
Y_{1,8}Y_{2,11}Y^{-1}_{2,13}Y^{-1}_{3,12}Y_{4,10}
$\\Monomial 188: $ (t^{-1} +t)$ $
Y^{-1}_{1,10}Y_{2,9}Y_{2,11}Y^{-1}_{2,13}Y^{-1}_{3,12}Y_{4,10}
$\\Monomial 193: $ (t^{-1} +t)$ $
Y_{1,8}Y_{2,11}Y^{-1}_{2,13}Y^{-1}_{4,14}
$\\Monomial 206: $ (t^{-1} +t)$ $
Y^{-1}_{1,10}Y_{2,9}Y_{2,11}Y^{-1}_{2,13}Y^{-1}_{4,14}$

Fundamental representation 3 :

\noindent Monomial 64: $ (t^{-1} +t)$ $
Y_{1,3}Y_{1,9}Y_{2,6}Y^{-1}_{2,8}Y^{-1}_{2,10}Y_{4,5}
$\\Monomial 90: $ (t^{-1} +t)$ $
Y^{-1}_{1,5}Y_{1,9}Y_{2,4}Y_{2,6}Y^{-1}_{2,8}Y^{-1}_{2,10}Y_{4,5}
$\\Monomial 91: $ (t^{-1} +t)$ $
Y_{3,5}Y^{-1}_{3,9}Y_{4,5}
$\\Monomial 93: $ (t^{-1} +t)$ $
Y_{1,3}Y^{-1}_{1,11}Y_{2,6}Y^{-1}_{2,8}Y_{4,5}
$\\Monomial 96: $ (t^{-1} +t)$ $
Y_{1,3}Y_{1,9}Y_{2,6}Y^{-1}_{2,8}Y^{-1}_{2,10}Y_{3,7}Y^{-1}_{4,9}
$\\Monomial 117: $ (t^{-1} +t)$ $
Y^{-1}_{1,5}Y^{-1}_{1,11}Y_{2,4}Y_{2,6}Y^{-1}_{2,8}Y_{4,5}
$\\Monomial 125: $ (t^{-1} +t)$ $
Y^{-1}_{1,5}Y_{1,9}Y_{2,4}Y_{2,6}Y^{-1}_{2,8}Y^{-1}_{2,10}Y_{3,7}Y^{-1}_{4,9}
$\\Monomial 126: $ (t^{-1} +t)$ $
Y_{3,5}Y_{3,7}Y^{-1}_{3,9}Y^{-1}_{4,9}
$\\Monomial 128: $ (t^{-1} +t)$ $
Y_{1,3}Y^{-1}_{1,11}Y_{2,6}Y^{-1}_{2,8}Y_{3,7}Y^{-1}_{4,9}
$\\Monomial 138: $ (t^{-1} +t)$ $
Y_{1,3}Y_{1,9}Y_{2,6}Y^{-1}_{3,11}
$\\Monomial 152: $ (t^{-1} +t)$ $
Y^{-1}_{1,5}Y^{-1}_{1,11}Y_{2,4}Y_{2,6}Y^{-1}_{2,8}Y_{3,7}Y^{-1}_{4,9}
$\\Monomial 159: $ (t^{-1} +t)$ $
Y_{2,8}Y_{2,10}Y_{3,5}Y^{-1}_{3,9}Y^{-1}_{3,11}
$\\Monomial 162: $ (t^{-1} +t)$ $
Y_{1,3}Y_{1,7}Y_{1,9}Y^{-1}_{2,8}Y_{3,7}Y^{-1}_{3,11}
$\\Monomial 165: $ (t^{-1} +t)$ $
Y_{1,3}Y^{-1}_{1,11}Y_{2,6}Y_{2,10}Y^{-1}_{3,11}
$\\Monomial 166: $ (t^{-1} +t)$ $
Y^{-1}_{1,5}Y_{1,9}Y_{2,4}Y_{2,6}Y^{-1}_{3,11}
$\\Monomial 194: $ (t^{-1} +t)$ $
Y^{-1}_{1,5}Y^{-1}_{1,11}Y_{2,4}Y_{2,6}Y_{2,10}Y^{-1}_{3,11}
$\\Monomial 208: $ (t^{-1} +t)$ $
Y^{-1}_{1,5}Y_{1,7}Y_{1,9}Y_{2,4}Y^{-1}_{2,8}Y_{3,7}Y^{-1}_{3,11}
$\\Monomial 209: $ (t^{-1} +t)$ $
Y_{1,11}Y_{2,8}Y^{-1}_{2,12}Y_{3,5}Y^{-1}_{3,9}
$\\Monomial 220: $ (t^{-1} +t)$ $
Y_{1,3}Y_{2,6}Y^{-1}_{2,12}
$\\Monomial 221: $ (t^{-1} +t)$ $
Y_{1,3}Y_{1,7}Y^{-1}_{1,11}Y^{-1}_{2,8}Y_{2,10}Y_{3,7}Y^{-1}_{3,11}
$\\Monomial 237: $ (t^{-1} +t)$ $
Y_{1,9}Y_{1,11}Y^{-1}_{2,10}Y^{-1}_{2,12}Y_{3,5}
$\\Monomial 238: $ (t^{-1} +t)$ $
Y_{1,7}Y_{1,9}Y^{-1}_{2,6}Y^{-1}_{2,8}Y_{3,5}Y_{3,7}Y^{-1}_{3,11}
$\\Monomial 251: $ (t^{-1} +t)$ $
Y^{-1}_{1,5}Y_{2,4}Y_{2,6}Y^{-1}_{2,12}
$\\Monomial 252: $ (t^{-1} +t)$ $
Y_{1,3}Y^{-1}_{1,9}Y^{-1}_{1,11}Y_{2,10}Y_{3,7}Y^{-1}_{3,11}
$\\Monomial 253: $ (t^{-1} +t)$ $
Y^{-1}_{1,5}Y_{1,7}Y^{-1}_{1,11}Y_{2,4}Y^{-1}_{2,8}Y_{2,10}Y_{3,7}Y^{-1}_{3,11}
$\\Monomial 257: $ (t^{-1} +t)$ $
Y_{1,3}Y_{1,7}Y^{-1}_{2,8}Y^{-1}_{2,12}Y_{3,7}
$\\Monomial 281: $ (t^{-1} +t)$ $
Y^{-1}_{1,13}Y_{2,8}Y_{3,5}Y^{-1}_{3,9}
$\\Monomial 289: $ (t^{-1} +t)$ $
Y^{-1}_{1,5}Y^{-1}_{1,9}Y^{-1}_{1,11}Y_{2,4}Y_{2,10}Y_{3,7}Y^{-1}_{3,11}
$\\Monomial 294: $ (t^{-1} +t)$ $
Y_{1,7}Y_{1,9}Y_{3,7}Y^{-1}_{3,9}Y^{-1}_{3,11}Y_{4,7}
$\\Monomial 296: $ (t^{-1} +t)$ $
Y_{1,9}Y_{1,11}Y_{2,6}Y_{2,8}Y^{-1}_{2,10}Y^{-1}_{2,12}Y^{-1}_{3,9}Y_{4,7}
$\\Monomial 298: $ (t^{-1} +t)$ $
Y_{1,3}Y^{-1}_{1,9}Y^{-1}_{2,12}Y_{3,7}
$\\Monomial 300: $ (t^{-1} +t)$ $
Y^{-1}_{1,5}Y_{1,7}Y_{2,4}Y^{-1}_{2,8}Y^{-1}_{2,12}Y_{3,7}
$\\Monomial 303: $ (t^{-1} +t)$ $
Y_{1,7}Y^{-1}_{1,11}Y^{-1}_{2,6}Y^{-1}_{2,8}Y_{2,10}Y_{3,5}Y_{3,7}Y^{-1}_{3,11}
$\\Monomial 320: $ (t^{-1} +t)$ $
Y_{1,9}Y^{-1}_{1,13}Y^{-1}_{2,10}Y_{3,5}
$\\Monomial 332: $ (t^{-1} +t)$ $
Y_{1,3}Y_{1,7}Y_{2,10}Y^{-1}_{2,12}Y^{-1}_{3,11}Y_{4,9}
$\\Monomial 351: $ (t^{-1} +t)$ $
Y^{-1}_{1,5}Y^{-1}_{1,9}Y_{2,4}Y^{-1}_{2,12}Y_{3,7}
$\\Monomial 353: $ (t^{-1} +t)$ $
Y_{1,7}Y^{-1}_{2,6}Y^{-1}_{2,8}Y^{-1}_{2,12}Y_{3,5}Y_{3,7}
$\\Monomial 359: $ (t^{-1} +t)$ $
Y^{-1}_{1,9}Y^{-1}_{1,11}Y^{-1}_{2,6}Y_{2,10}Y_{3,5}Y_{3,7}Y^{-1}_{3,11}
$\\Monomial 361: $ (t^{-1} +t)$ $
Y_{1,7}Y^{-1}_{1,11}Y_{2,10}Y_{3,7}Y^{-1}_{3,9}Y^{-1}_{3,11}Y_{4,7}
$\\Monomial 362: $ (t^{-1} +t)$ $
Y^{-1}_{1,11}Y^{-1}_{1,13}Y_{3,5}
$\\Monomial 368: $ (t^{-1} +t)$ $
Y_{1,7}Y_{1,9}Y_{3,7}Y^{-1}_{3,11}Y^{-1}_{4,11}
$\\Monomial 370: $ (t^{-1} +t)$ $
Y_{1,9}Y_{1,11}Y_{2,6}Y_{2,8}Y^{-1}_{2,10}Y^{-1}_{2,12}Y^{-1}_{4,11}
$\\Monomial 382: $ (t^{-1} +t)$ $
Y_{1,3}Y^{-1}_{1,9}Y_{2,8}Y_{2,10}Y^{-1}_{2,12}Y^{-1}_{3,11}Y_{4,9}
$\\Monomial 384: $ (t^{-1} +t)$ $
Y^{-1}_{1,5}Y_{1,7}Y_{2,4}Y_{2,10}Y^{-1}_{2,12}Y^{-1}_{3,11}Y_{4,9}
$\\Monomial 394: $ (t^{-1} +t)$ $
Y_{1,9}Y^{-1}_{1,13}Y_{2,6}Y_{2,8}Y^{-1}_{2,10}Y^{-1}_{3,9}Y_{4,7}
$\\Monomial 399: $ (t^{-1} +t)$ $
Y_{1,3}Y_{1,7}Y_{2,10}Y^{-1}_{2,12}Y^{-1}_{4,13}
$\\Monomial 414: $ (t^{-1} +t)$ $
Y^{-1}_{1,9}Y^{-1}_{2,6}Y^{-1}_{2,12}Y_{3,5}Y_{3,7}
$\\Monomial 422: $ (t^{-1} +t)$ $
Y^{-1}_{1,5}Y^{-1}_{1,9}Y_{2,4}Y_{2,8}Y_{2,10}Y^{-1}_{2,12}Y^{-1}_{3,11}Y_{4,9}
$\\Monomial 428: $ (t^{-1} +t)$ $
Y^{-1}_{1,9}Y^{-1}_{1,11}Y_{2,8}Y_{2,10}Y_{3,7}Y^{-1}_{3,9}Y^{-1}_{3,11}Y_{4,7}
$\\Monomial 431: $ (t^{-1} +t)$ $
Y_{1,7}Y^{-1}_{1,11}Y_{2,10}Y_{3,7}Y^{-1}_{3,11}Y^{-1}_{4,11}
$\\Monomial 432: $ (t^{-1} +t)$ $
Y^{-1}_{1,11}Y^{-1}_{1,13}Y_{2,6}Y_{2,8}Y^{-1}_{3,9}Y_{4,7}
$\\Monomial 436: $ (t^{-1} +t)$ $
Y_{1,7}Y^{-1}_{2,6}Y_{2,10}Y^{-1}_{2,12}Y_{3,5}Y^{-1}_{3,11}Y_{4,9}
$\\Monomial 438: $ (t^{-1} +t)$ $
Y_{1,7}Y^{-1}_{2,12}Y_{3,7}Y^{-1}_{3,9}Y_{4,7}
$\\Monomial 461: $ (t^{-1} +t)$ $
Y_{1,3}Y^{-1}_{1,9}Y_{2,8}Y_{2,10}Y^{-1}_{2,12}Y^{-1}_{4,13}
$\\Monomial 463: $ (t^{-1} +t)$ $
Y^{-1}_{1,5}Y_{1,7}Y_{2,4}Y_{2,10}Y^{-1}_{2,12}Y^{-1}_{4,13}
$\\Monomial 469: $ (t^{-1} +t)$ $
Y_{1,9}Y^{-1}_{1,13}Y_{2,6}Y_{2,8}Y^{-1}_{2,10}Y^{-1}_{4,11}
$\\Monomial 495: $ (t^{-1} +t)$ $
Y^{-1}_{1,5}Y^{-1}_{1,9}Y_{2,4}Y_{2,8}Y_{2,10}Y^{-1}_{2,12}Y^{-1}_{4,13}
$\\Monomial 498: $ (t^{-1} +t)$ $
Y_{1,9}Y^{-1}_{1,11}Y^{-1}_{1,13}Y_{2,6}Y^{-1}_{2,10}Y_{4,7}
$\\Monomial 500: $ (t^{-1} +t)$ $
Y^{-1}_{1,9}Y_{2,8}Y^{-1}_{2,12}Y_{3,7}Y^{-1}_{3,9}Y_{4,7}
$\\Monomial 502: $ (t^{-1} +t)$ $
Y_{1,7}Y_{2,8}Y_{2,10}Y^{-1}_{2,12}Y^{-1}_{3,9}Y^{-1}_{3,11}Y_{4,7}Y_{4,9}
$\\Monomial 505: $ (t^{-1} +t)$ $
Y^{-1}_{1,9}Y^{-1}_{2,6}Y_{2,8}Y_{2,10}Y^{-1}_{2,12}Y_{3,5}Y^{-1}_{3,11}Y_{4,9}
$\\Monomial 512: $ (t^{-1} +t)$ $
Y^{-1}_{1,9}Y^{-1}_{1,11}Y_{2,8}Y_{2,10}Y_{3,7}Y^{-1}_{3,11}Y^{-1}_{4,11}
$\\Monomial 514: $ (t^{-1} +t)$ $
Y^{-1}_{1,11}Y^{-1}_{1,13}Y_{2,6}Y_{2,8}Y^{-1}_{4,11}
$\\Monomial 526: $ (t^{-1} +t)$ $
Y_{1,7}Y^{-1}_{2,6}Y_{2,10}Y^{-1}_{2,12}Y_{3,5}Y^{-1}_{4,13}
$\\Monomial 528: $ (t^{-1} +t)$ $
Y_{1,7}Y^{-1}_{2,12}Y_{3,7}Y^{-1}_{4,11}
$\\Monomial 564: $ (t^{-1} +t)$ $
Y^{-1}_{2,10}Y^{-1}_{2,12}Y_{3,7}Y_{4,7}
$\\Monomial 575: $ (t^{-1} +t)$ $
Y_{1,7}Y_{1,9}Y^{-1}_{1,11}Y^{-1}_{1,13}Y^{-1}_{2,8}Y^{-1}_{2,10}Y_{3,7}Y_{4,7}
$\\Monomial 577: $ (t^{-1} +t)$ $
Y^{-1}_{1,9}Y^2_{2,8}Y_{2,10}Y^{-1}_{2,12}Y^{-1}_{3,9}Y^{-1}_{3,11}Y_{4,7}Y_{4,9}
$\\Monomial 581: $ (t^{-1} +t)$ $
Y_{1,9}Y^{-1}_{1,11}Y^{-1}_{1,13}Y_{2,6}Y^{-1}_{2,10}Y_{3,9}Y^{-1}_{4,11}
$\\Monomial 583: $ (t^{-1} +t)$ $
Y^{-1}_{1,9}Y_{2,8}Y^{-1}_{2,12}Y_{3,7}Y^{-1}_{4,11}
$\\Monomial 586: $ (t^{-1} +t)$ $
Y_{1,7}Y_{2,8}Y_{2,10}Y^{-1}_{2,12}Y^{-1}_{3,9}Y_{4,7}Y^{-1}_{4,13}
$\\Monomial 591: $ (t^{-1} +t)$ $
Y^{-1}_{1,9}Y^{-1}_{2,6}Y_{2,8}Y_{2,10}Y^{-1}_{2,12}Y_{3,5}Y^{-1}_{4,13}
$\\Monomial 622: $ (t^{-1} +t)$ $
Y_{1,7}Y_{2,8}Y_{2,10}Y^{-1}_{2,12}Y^{-1}_{3,11}Y_{4,9}Y^{-1}_{4,11}
$\\Monomial 648: $ (t^{-1} +t)$ $
Y^{-1}_{2,10}Y^{-1}_{2,12}Y_{3,7}Y_{3,9}Y^{-1}_{4,11}
$\\Monomial 656: $ (t^{-1} +t)$ $
Y_{2,8}Y^{-1}_{2,12}Y^{-1}_{3,11}Y_{4,7}Y_{4,9}
$\\Monomial 657: $ (t^{-1} +t)$ $
Y_{1,7}Y_{1,9}Y^{-1}_{1,11}Y^{-1}_{1,13}Y^{-1}_{2,8}Y^{-1}_{2,10}Y_{3,7}Y_{3,9}Y^{-1}_{4,11}
$\\Monomial 666: $ (t^{-1} +t)$ $
Y_{1,7}Y_{1,9}Y^{-1}_{1,11}Y^{-1}_{1,13}Y^{-1}_{3,11}Y_{4,7}Y_{4,9}
$\\Monomial 669: $ (t^{-1} +t)$ $
Y_{1,7}Y_{2,8}Y_{2,10}Y^{-1}_{2,12}Y^{-1}_{4,11}Y^{-1}_{4,13}
$\\Monomial 672: $ (t^{-1} +t)$ $
Y^{-1}_{1,9}Y^2_{2,8}Y_{2,10}Y^{-1}_{2,12}Y^{-1}_{3,9}Y_{4,7}Y^{-1}_{4,13}
$\\Monomial 694: $ (t^{-1} +t)$ $
Y_{1,9}Y^{-1}_{1,11}Y^{-1}_{1,13}Y_{2,6}Y_{2,12}Y^{-1}_{3,13}
$\\Monomial 696: $ (t^{-1} +t)$ $
Y^{-1}_{1,9}Y^2_{2,8}Y_{2,10}Y^{-1}_{2,12}Y^{-1}_{3,11}Y_{4,9}Y^{-1}_{4,11}
$\\Monomial 729: $ (t^{-1} +t)$ $
Y_{1,11}Y_{2,8}Y^{-1}_{2,10}Y^{-2}_{2,12}Y_{4,7}Y_{4,9}
$\\Monomial 755: $ (t^{-1} +t)$ $
Y^{-1}_{1,9}Y^2_{2,8}Y_{2,10}Y^{-1}_{2,12}Y^{-1}_{4,11}Y^{-1}_{4,13}
$\\Monomial 764: $ (t^{-1} +t)$ $
Y_{2,8}Y^{-1}_{2,12}Y_{4,7}Y^{-1}_{4,13}
$\\Monomial 765: $ (t^{-1} +t)$ $
Y_{1,7}Y_{1,9}Y^{-1}_{1,11}Y^{-1}_{1,13}Y_{4,7}Y^{-1}_{4,13}
$\\Monomial 767: $ (t^{-1} +t)$ $
Y_{2,8}Y^{-1}_{2,12}Y_{3,9}Y^{-1}_{3,11}Y_{4,9}Y^{-1}_{4,11}
$\\Monomial 768: $ (t^{-1} +t)$ $
Y_{1,7}Y_{1,9}Y^{-1}_{1,11}Y^{-1}_{1,13}Y_{3,9}Y^{-1}_{3,11}Y_{4,9}Y^{-1}_{4,11}
Y_{2,6}Y_{2,8}Y^{-1}_{2,12}Y^{-1}_{2,14}Y^{-1}_{3,9}Y_{3,11}
$\\Monomial 770: $ (t^{-1} +t)$ $
Y_{1,9}Y^{-1}_{1,11}Y_{2,6}Y^{-1}_{2,14}
$\\Monomial 771: $ (t^{-1} +t)$ $
Y_{1,7}Y_{1,9}Y^{-1}_{1,11}Y^{-1}_{1,13}Y^{-1}_{2,8}Y_{2,12}Y_{3,7}Y^{-1}_{3,13}
$\\Monomial 772: $ (t^{-1} +t)$ $
Y_{3,7}Y^{-1}_{3,13}
$\\Monomial 815: $ (t^{-1} +t)$ $
Y_{1,11}Y_{2,8}Y^{-1}_{2,10}Y^{-2}_{2,12}Y_{3,11}Y_{4,7}Y^{-1}_{4,13}
$\\Monomial 818: $ (t^{-1} +t)$ $
Y_{1,11}Y_{2,8}Y^{-1}_{2,10}Y^{-2}_{2,12}Y_{3,9}Y_{4,9}Y^{-1}_{4,11}
$\\Monomial 822: $ (t^{-1} +t)$ $
Y_{1,7}Y_{1,9}Y^{-1}_{1,11}Y^{-1}_{2,8}Y^{-1}_{2,14}Y_{3,7}
$\\Monomial 834: $ (t^{-1} +t)$ $
Y^{-1}_{1,13}Y_{2,8}Y^{-1}_{2,10}Y^{-1}_{2,12}Y_{4,7}Y_{4,9}
$\\Monomial 840: $ (t^{-1} +t)$ $
Y_{2,8}Y_{2,10}Y^{-1}_{3,11}Y^{-1}_{3,13}Y_{4,9}
$\\Monomial 844: $ (t^{-1} +t)$ $
Y_{2,8}Y^{-1}_{2,12}Y_{3,9}Y^{-1}_{4,11}Y^{-1}_{4,13}
$\\Monomial 845: $ (t^{-1} +t)$ $
Y_{1,7}Y_{1,9}Y^{-1}_{1,11}Y^{-1}_{1,13}Y_{2,10}Y_{2,12}Y^{-1}_{3,11}Y^{-1}_{3,13}Y_{4,9}
$\\Monomial 849: $ (t^{-1} +t)$ $
Y_{1,7}Y_{1,9}Y^{-1}_{1,11}Y^{-1}_{1,13}Y_{3,9}Y^{-1}_{4,11}Y^{-1}_{4,13}
$\\Monomial 907: $ (t^{-1} +t)$ $
Y_{1,11}Y_{2,8}Y^{-1}_{2,10}Y^{-2}_{2,12}Y_{3,9}Y_{3,11}Y^{-1}_{4,11}Y^{-1}_{4,13}
$\\Monomial 911: $ (t^{-1} +t)$ $
Y_{1,11}Y_{2,8}Y^{-1}_{2,10}Y^{-1}_{2,12}Y_{2,14}Y^{-1}_{3,15}Y_{4,7}
$\\Monomial 916: $ (t^{-1} +t)$ $
Y_{2,8}Y_{2,10}Y^{-1}_{3,13}Y^{-1}_{4,13}
$\\Monomial 920: $ (t^{-1} +t)$ $
Y_{1,11}Y_{2,8}Y^{-1}_{2,12}Y^{-1}_{3,13}Y_{4,9}
$\\Monomial 928: $ (t^{-1} +t)$ $
Y^{-1}_{1,13}Y_{2,8}Y^{-1}_{2,10}Y^{-1}_{2,12}Y_{3,11}Y_{4,7}Y^{-1}_{4,13}
$\\Monomial 930: $ (t^{-1} +t)$ $
Y^{-1}_{1,13}Y_{2,8}Y^{-1}_{2,10}Y^{-1}_{2,12}Y_{3,9}Y_{4,9}Y^{-1}_{4,11}
$\\Monomial 950: $ (t^{-1} +t)$ $
Y_{1,7}Y_{1,9}Y^{-1}_{1,11}Y^{-1}_{1,13}Y_{2,10}Y_{2,12}Y^{-1}_{3,13}Y^{-1}_{4,13}
$\\Monomial 953: $ (t^{-1} +t)$ $
Y_{1,7}Y_{1,9}Y^{-1}_{1,11}Y_{2,10}Y^{-1}_{2,14}Y^{-1}_{3,11}Y_{4,9}
$\\Monomial 979: $ (t^{-1} +t)$ $
Y_{1,11}Y_{1,15}Y_{2,8}Y^{-1}_{2,10}Y^{-1}_{2,12}Y^{-1}_{2,16}Y_{4,7}
$\\Monomial 981: $ (t^{-1} +t)$ $
Y_{1,9}Y_{1,11}Y^{-1}_{2,10}Y^{-1}_{2,12}Y_{3,9}Y^{-1}_{3,13}Y_{4,9}
$\\Monomial 998: $ (t^{-1} +t)$ $
Y^{-1}_{1,13}Y_{2,8}Y^{-1}_{2,10}Y_{2,14}Y^{-1}_{3,15}Y_{4,7}
$\\Monomial 1001: $ (t^{-1} +t)$ $
Y^{-1}_{1,13}Y_{2,8}Y^{-1}_{2,10}Y^{-1}_{2,12}Y_{3,9}Y_{3,11}Y^{-1}_{4,11}Y^{-1}_{4,13}
$\\Monomial 1005: $ (t^{-1} +t)$ $
Y_{1,11}Y_{2,8}Y^{-1}_{2,12}Y_{3,11}Y^{-1}_{3,13}Y^{-1}_{4,13}
$\\Monomial 1016: $ (t^{-1} +t)$ $
Y_{1,11}Y_{2,8}Y^{-1}_{2,10}Y^{-1}_{2,12}Y_{2,14}Y_{3,9}Y^{-1}_{3,15}Y^{-1}_{4,11}
$\\Monomial 1017: $ (t^{-1} +t)$ $
Y_{1,7}Y_{1,9}Y^{-1}_{2,12}Y^{-1}_{2,14}Y_{4,9}
$\\Monomial 1026: $ (t^{-1} +t)$ $
Y^{-1}_{1,13}Y_{2,8}Y^{-1}_{3,13}Y_{4,9}
$\\Monomial 1044: $ (t^{-1} +t)$ $
Y_{1,7}Y_{1,9}Y^{-1}_{1,11}Y_{2,10}Y^{-1}_{2,14}Y^{-1}_{4,13}
$\\Monomial 1058: $ (t^{-1} +t)$ $
Y_{1,11}Y_{1,15}Y_{2,8}Y^{-1}_{2,10}Y^{-1}_{2,12}Y^{-1}_{2,16}Y_{3,9}Y^{-1}_{4,11}
$\\Monomial 1060: $ (t^{-1} +t)$ $
Y_{1,9}Y_{1,11}Y^{-1}_{2,10}Y^{-1}_{2,12}Y_{3,9}Y_{3,11}Y^{-1}_{3,13}Y^{-1}_{4,13}
$\\Monomial 1074: $ (t^{-1} +t)$ $
Y_{1,11}Y^{-1}_{1,17}Y_{2,8}Y^{-1}_{2,10}Y^{-1}_{2,12}Y_{4,7}
$\\Monomial 1075: $ (t^{-1} +t)$ $
Y^{-1}_{1,13}Y_{1,15}Y_{2,8}Y^{-1}_{2,10}Y^{-1}_{2,16}Y_{4,7}
$\\Monomial 1076: $ (t^{-1} +t)$ $
Y_{1,7}Y^{-1}_{1,11}Y_{2,10}Y^{-1}_{2,12}Y^{-1}_{2,14}Y_{4,9}
$\\Monomial 1078: $ (t^{-1} +t)$ $
Y_{1,9}Y^{-1}_{1,13}Y^{-1}_{2,10}Y_{3,9}Y^{-1}_{3,13}Y_{4,9}
$\\Monomial 1079: $ (t^{-1} +t)$ $
Y_{1,11}Y_{2,8}Y_{2,14}Y^{-1}_{3,13}Y^{-1}_{3,15}
$\\Monomial 1085: $ (t^{-1} +t)$ $
Y^{-1}_{1,13}Y_{2,8}Y_{3,11}Y^{-1}_{3,13}Y^{-1}_{4,13}
$\\Monomial 1096: $ (t^{-1} +t)$ $
Y^{-1}_{1,13}Y_{2,8}Y^{-1}_{2,10}Y_{2,14}Y_{3,9}Y^{-1}_{3,15}Y^{-1}_{4,11}
$\\Monomial 1101: $ (t^{-1} +t)$ $
Y_{1,7}Y_{1,9}Y^{-1}_{2,12}Y^{-1}_{2,14}Y_{3,11}Y^{-1}_{4,13}
$\\Monomial 1123: $ (t^{-1} +t)$ $
Y^{-1}_{1,13}Y^{-1}_{1,17}Y_{2,8}Y^{-1}_{2,10}Y_{4,7}
$\\Monomial 1137: $ (t^{-1} +t)$ $
Y^{-1}_{1,9}Y^{-1}_{1,11}Y_{2,8}Y_{2,10}Y^{-1}_{2,12}Y^{-1}_{2,14}Y_{4,9}
$\\Monomial 1138: $ (t^{-1} +t)$ $
Y^{-1}_{1,11}Y^{-1}_{1,13}Y_{3,9}Y^{-1}_{3,13}Y_{4,9}
$\\Monomial 1140: $ (t^{-1} +t)$ $
Y_{1,9}Y_{1,11}Y^{-1}_{2,10}Y_{2,14}Y_{3,9}Y^{-1}_{3,13}Y^{-1}_{3,15}
$\\Monomial 1146: $ (t^{-1} +t)$ $
Y_{1,11}Y_{1,15}Y_{2,8}Y^{-1}_{2,16}Y^{-1}_{3,13}
$\\Monomial 1154: $ (t^{-1} +t)$ $
Y_{1,11}Y^{-1}_{1,17}Y_{2,8}Y^{-1}_{2,10}Y^{-1}_{2,12}Y_{3,9}Y^{-1}_{4,11}
$\\Monomial 1155: $ (t^{-1} +t)$ $
Y^{-1}_{1,13}Y_{1,15}Y_{2,8}Y^{-1}_{2,10}Y^{-1}_{2,16}Y_{3,9}Y^{-1}_{4,11}
$\\Monomial 1156: $ (t^{-1} +t)$ $
Y_{1,7}Y^{-1}_{1,11}Y_{2,10}Y^{-1}_{2,12}Y^{-1}_{2,14}Y_{3,11}Y^{-1}_{4,13}
$\\Monomial 1158: $ (t^{-1} +t)$ $
Y_{1,9}Y^{-1}_{1,13}Y^{-1}_{2,10}Y_{3,9}Y_{3,11}Y^{-1}_{3,13}Y^{-1}_{4,13}
$\\Monomial 1160: $ (t^{-1} +t)$ $
Y^{-1}_{1,13}Y_{2,8}Y_{2,12}Y_{2,14}Y^{-1}_{3,13}Y^{-1}_{3,15}
$\\Monomial 1177: $ (t^{-1} +t)$ $
Y_{1,7}Y_{1,9}Y^{-1}_{3,15}
$\\Monomial 1191: $ (t^{-1} +t)$ $
Y_{1,9}Y_{1,11}Y_{1,15}Y^{-1}_{2,10}Y^{-1}_{2,16}Y_{3,9}Y^{-1}_{3,13}
$\\Monomial 1193: $ (t^{-1} +t)$ $
Y^{-1}_{1,13}Y^{-1}_{1,17}Y_{2,8}Y^{-1}_{2,10}Y_{3,9}Y^{-1}_{4,11}
$\\Monomial 1204: $ (t^{-1} +t)$ $
Y^{-1}_{1,9}Y^{-1}_{1,11}Y_{2,8}Y_{2,10}Y^{-1}_{2,12}Y^{-1}_{2,14}Y_{3,11}Y^{-1}_{4,13}
$\\Monomial 1205: $ (t^{-1} +t)$ $
Y^{-1}_{1,11}Y^{-1}_{1,13}Y_{3,9}Y_{3,11}Y^{-1}_{3,13}Y^{-1}_{4,13}
$\\Monomial 1209: $ (t^{-1} +t)$ $
Y_{1,9}Y^{-1}_{1,13}Y^{-1}_{2,10}Y_{2,12}Y_{2,14}Y_{3,9}Y^{-1}_{3,13}Y^{-1}_{3,15}
$\\Monomial 1231: $ (t^{-1} +t)$ $
Y_{1,11}Y^{-1}_{1,17}Y_{2,8}Y^{-1}_{3,13}
$\\Monomial 1232: $ (t^{-1} +t)$ $
Y^{-1}_{1,13}Y_{1,15}Y_{2,8}Y_{2,12}Y^{-1}_{2,16}Y^{-1}_{3,13}
$\\Monomial 1239: $ (t^{-1} +t)$ $
Y_{1,7}Y^{-1}_{1,11}Y_{2,10}Y^{-1}_{3,15}
$\\Monomial 1256: $ (t^{-1} +t)$ $
Y^{-1}_{1,13}Y^{-1}_{1,17}Y_{2,8}Y_{2,12}Y^{-1}_{3,13}
$\\Monomial 1263: $ (t^{-1} +t)$ $
Y^{-1}_{1,11}Y^{-1}_{1,13}Y_{2,12}Y_{2,14}Y_{3,9}Y^{-1}_{3,13}Y^{-1}_{3,15}
$\\Monomial 1276: $ (t^{-1} +t)$ $
Y_{1,7}Y^{-1}_{2,12}Y_{3,11}Y^{-1}_{3,15}
$\\Monomial 1277: $ (t^{-1} +t)$ $
Y_{1,15}Y_{2,8}Y^{-1}_{2,14}Y^{-1}_{2,16}
$\\Monomial 1278: $ (t^{-1} +t)$ $
Y_{1,9}Y^{-1}_{1,13}Y_{1,15}Y^{-1}_{2,10}Y_{2,12}Y^{-1}_{2,16}Y_{3,9}Y^{-1}_{3,13}
$\\Monomial 1284: $ (t^{-1} +t)$ $
Y_{1,9}Y_{1,11}Y^{-1}_{1,17}Y^{-1}_{2,10}Y_{3,9}Y^{-1}_{3,13}
$\\Monomial 1288: $ (t^{-1} +t)$ $
Y^{-1}_{1,9}Y^{-1}_{1,11}Y_{2,8}Y_{2,10}Y^{-1}_{3,15}
$\\Monomial 1305: $ (t^{-1} +t)$ $
Y_{1,9}Y_{1,15}Y^{-1}_{2,10}Y^{-1}_{2,14}Y^{-1}_{2,16}Y_{3,9}
$\\Monomial 1310: $ (t^{-1} +t)$ $
Y^{-1}_{1,17}Y_{2,8}Y^{-1}_{2,14}
$\\Monomial 1311: $ (t^{-1} +t)$ $
Y_{1,9}Y^{-1}_{1,13}Y^{-1}_{1,17}Y^{-1}_{2,10}Y_{2,12}Y_{3,9}Y^{-1}_{3,13}
$\\Monomial 1312: $ (t^{-1} +t)$ $
Y^{-1}_{1,11}Y^{-1}_{1,13}Y_{1,15}Y_{2,12}Y^{-1}_{2,16}Y_{3,9}Y^{-1}_{3,13}
$\\Monomial 1317: $ (t^{-1} +t)$ $
Y^{-1}_{1,9}Y_{2,8}Y^{-1}_{2,12}Y_{3,11}Y^{-1}_{3,15}
$\\Monomial 1346: $ (t^{-1} +t)$ $
Y^{-1}_{1,11}Y^{-1}_{1,13}Y^{-1}_{1,17}Y_{2,12}Y_{3,9}Y^{-1}_{3,13}
$\\Monomial 1348: $ (t^{-1} +t)$ $
Y^{-1}_{1,11}Y_{1,15}Y^{-1}_{2,14}Y^{-1}_{2,16}Y_{3,9}
$\\Monomial 1349: $ (t^{-1} +t)$ $
Y^{-1}_{1,10}Y^{-1}_{1,12}Y_{2,9}Y_{2,11}Y^{-1}_{2,15}
$\\Monomial 1356: $ (t^{-1} +t)$ $
Y_{1,9}Y^{-1}_{1,17}Y^{-1}_{2,10}Y^{-1}_{2,14}Y_{3,9}
$\\Monomial 1360: $ (t^{-1} +t)$ $
Y_{1,9}Y_{1,15}Y_{2,12}Y^{-1}_{2,14}Y^{-1}_{2,16}Y^{-1}_{3,13}Y_{4,11}
$\\Monomial 1387: $ (t^{-1} +t)$ $
Y^{-1}_{1,11}Y^{-1}_{1,17}Y^{-1}_{2,14}Y_{3,9}
$\\Monomial 1392: $ (t^{-1} +t)$ $
Y_{3,11}Y^{-1}_{3,13}Y^{-1}_{3,15}Y_{4,11}
$\\Monomial 1394: $ (t^{-1} +t)$ $
Y^{-1}_{1,11}Y_{1,15}Y_{2,10}Y_{2,12}Y^{-1}_{2,14}Y^{-1}_{2,16}Y^{-1}_{3,13}Y_{4,11}
$\\Monomial 1397: $ (t^{-1} +t)$ $
Y_{1,9}Y^{-1}_{1,17}Y_{2,12}Y^{-1}_{2,14}Y^{-1}_{3,13}Y_{4,11}
$\\Monomial 1398: $ (t^{-1} +t)$ $
Y_{1,9}Y_{1,15}Y_{2,12}Y^{-1}_{2,14}Y^{-1}_{2,16}Y^{-1}_{4,15}
$\\Monomial 1424: $ (t^{-1} +t)$ $
Y^{-1}_{1,11}Y^{-1}_{1,17}Y_{2,10}Y_{2,12}Y^{-1}_{2,14}Y^{-1}_{3,13}Y_{4,11}
$\\Monomial 1431: $ (t^{-1} +t)$ $
Y_{3,11}Y^{-1}_{3,15}Y^{-1}_{4,15}
$\\Monomial 1433: $ (t^{-1} +t)$ $
Y^{-1}_{1,11}Y_{1,15}Y_{2,10}Y_{2,12}Y^{-1}_{2,14}Y^{-1}_{2,16}Y^{-1}_{4,15}
$\\Monomial 1436: $ (t^{-1} +t)$ $
Y_{1,9}Y^{-1}_{1,17}Y_{2,12}Y^{-1}_{2,14}Y^{-1}_{4,15}
$\\Monomial 1452: $ (t^{-1} +t)$ $
Y^{-1}_{1,11}Y^{-1}_{1,17}Y_{2,10}Y_{2,12}Y^{-1}_{2,14}Y^{-1}_{4,15}$

\end{document}